\documentclass[11pt]{article}
\usepackage[total={6.5in,9in},top=1in, left=1in, includefoot]{geometry}
\usepackage{amssymb,amsmath,amsthm,amscd,amsfonts}
\usepackage{mathptmx,microtype,booktabs}
\usepackage{graphicx,subfigure}
\usepackage{algorithmic}
\usepackage{multirow}
\usepackage{color}
\usepackage{tikz}
\usepackage{subfigure}
\numberwithin{equation}{section}

\newtheorem{theorem}{Theorem}[section]

\newtheorem{remark}[theorem]{Remark}

\newtheorem{algori}[theorem]{Algorithm}

\newcommand{\abs}[1]{\left|#1\right|}

\newcommand{\R}{\mathbb{R}}

\newcommand {\eps}  {\varepsilon}
\newcommand{\ds}{\displaystyle}
\newcommand{\ud}{\, \mathrm{d}}

\title{An Asymptotic Preserving Scheme for the Diffusive Limit \\ of Kinetic systems for Chemotaxis}

\author{Jos\'{e} A. Carrillo$^1$, Bokai Yan$^2$\\
\vspace{5pt}\small{$^{1}$ Instituci\'o Catalana de Recerca i Estudis Avan\c cats and Departament de Matem\`atiques}\\[-8pt]
\small{Universitat Aut\`onoma de Barcelona, E-08193 Bellaterra, Spain}\\[-4pt]
\small{Email: \texttt{carrillo@mat.uab.es}} \\
\small{{\it On leave from:} Department of Mathematics, Imperial
College London,
London SW7 2AZ, UK.}\\
\vspace{5pt}\small{$^{2}$ Department of Mathematics}\\[-8pt]
\small{University of Wisconsin-Madison, WI 53706}\\[-4pt]
\small{Email: \texttt{yan@math.wisc.edu}}}

\begin{document}

\maketitle

\begin{abstract}
In this work we numerically study the diffusive limit of run \&
tumble kinetic models for cell motion due to chemotaxis by means
of asymptotic preserving schemes. It is well-known that the
diffusive limit of these models leads to the classical
Patlak-Keller-Segel macroscopic model for chemotaxis. We will show
that the proposed scheme is able to accurately approximate the
solutions before blow-up time for small parameter. Moreover, the
numerical results indicate that the global solutions of the
kinetic models stabilize for long times to steady states for all
the analyzed parameter range. We also generalize these asymptotic
preserving schemes to two dimensional kinetic models in the radial
case. The blow-up of solutions is numerically investigated in all
these cases.
\end{abstract}

\section{Introduction}

Chemotaxis is one of the basic mechanisms of cell motility due to
chemical interaction. Cells are attracted by the concentration of
certain chemical substance, called chemoattractant, and they
direct their movement toward the regions of highest concentration.
Typically, this phenomenon have been described based on
macroscopic systems of equations describing the evolution of the
cell density and the chemoattractant in time. These systems of
drift-diffusion type are the well-know classical
(Patlak)-Keller-Segel models \cite{Pa53,KS70,KS71}. Another point
of view was introduced by a mesoscopic description of these
phenomena bridging from stochastic interacting particle systems to
macroscopic equations. This middle ground consists in describing
the movement of cells by a ``run \& tumble'' process
\cite{ODA88,OS97}. The cells move along a straight line in the
running phase and make reorientation as a reaction to the
surrounding chemicals during the tumbling phase. This is the
typical behavior that has been observed in experiments. The
resulting nondimensionalized kinetic equation, with parabolic
rescaling, reads
\begin{equation}
\label{eq:kinetic}
  \ds \eps \frac{{\partial f  }}{{\partial
t}} +  v  \cdot \nabla_x f     = \frac{1}{\eps}\int_V \left(
{T_\eps  f' - T_\eps ^* f} \right) \ud v'
 \end{equation}
Here $f=f_\eps(t,x,v)$ is the density of cells at position $x\in
\R^N$, moving with velocity $v\in V \subset \R^N$. $S=S_\eps(t,x)$
is the density of chemoattractant, whose governing equation will
be described later.  $T_\eps = T_\eps[S](t,x,v,v')$ is the turning
kernel operator. We use the abbreviation $f'=f(t,x,v')$,
$T_\eps^*[S](t,x,v,v') = T_\eps[S](t,x,v',v)$. $\eps$ is the ratio
of the mean running length between jumps to the typical
observation length scale. For simplicity we omit the dependence on
$\eps$ in notations (except $T_\eps$, where $\eps$ explicitly
enters into the expression.) We refer to \cite{CMPS04} for the
details on the rescaling. The chemoattractant $S$ is typically
given by the Poisson equation
\begin{equation}
\label{eq:S_Poisson} -\Delta S = \rho \quad \mbox{where }\, \rho =
\rho_\eps(t,x) = \ds\int_V f_\eps(t,x,v) \ud v \,\mbox{ is the
macroscopic density of cells.}
\end{equation}
Here, we assume that the time scale of relaxation for the
chemoattractant is much smaller that the one for the cell density.
Although this is the most standard way of introducing the coupling
with the chemoattractant, some authors have proposed to use some
generalized models for any dimension where the chemoattractant $S$
is given by
\begin{equation}
\label{eq:S_log} S = -\frac{1}{N\pi}\log|x|\ast\rho\,.
\end{equation}
Note that in two dimensional case, (\ref{eq:S_log}) is exactly
(\ref{eq:S_Poisson}).

Starting with the kinetic equation (\ref{eq:kinetic}), one can (at
lease formally) derive the macroscopic limit as $\eps\to0$,
\begin{equation}
  \label{eq:KellerSegel}
  \partial_t\rho = \nabla_x \cdot \left(D\nabla_x\rho - \chi\rho\nabla_x S\right).
\end{equation}
The details will be described in Section \ref{sec:macrolimit}.
Coupling with (\ref{eq:S_Poisson}), one obtains the already
mentioned Patlak-Keller-Segel system \cite{Pa53, KS70, KS71} for
chemotaxis. We refer to the reviews \cite{Hor02,Perbook} and the
references therein for mathematical results on this system. The
behavior of the solution is quite different depending on the
dimension. In the $1$D case, the global solution exists for all
initial condition. In the $2$D case, there exists a critical mass
$M_c$ determined by the coefficients \cite{DP05,BDP07}. The global
solution exists if the initial total mass satisfies $M<M_c$
(subcritical case) and its long time asymptotics are given by
self-similar profiles, otherwise the solution blows up in infinite
(critical case $M=M_c$ \cite{BCM08}) or finite time (supercritical
case $M>M_c$ \cite{DS09}). The critical case has infinitely many
stationary states, each of them having their own basin of
attraction \cite{BCC10} depending on the tail of the distribution
at infinity. In the $3$D case, the relevant quantity ensuring
global existence is the $L^{3/2}$ norm of the initial density, see
\cite{CPZ04}, together with some condition involving the second
moment of the solution.

Similarly, starting with (\ref{eq:kinetic})-(\ref{eq:S_log}), one
derives the modified Keller-Segel model
(\ref{eq:KellerSegel})-(\ref{eq:S_log}) studied in \cite{CPS07},
as $\eps\to0$. For the modified system, the existence of critical
mass is extended to the $1$D and $3$D cases. It is given by
\begin{equation*}
  M_c = \frac{2N^2\pi D}{\chi}.
\end{equation*}
A very interesting consequence is that even in one dimension, this
system blows-up over this critical mass. This makes possible to
numerically analyse this highly subtle phenomena. This modified
Keller-Segel system was numerically solved in \cite{BCC08} by
using optimal transportation ideas and the scheme proven to be
convergent in the subcritical case.

Alt in \cite{Alt80, Alt81} derived (\ref{eq:KellerSegel}) from a
transport equation similar to (\ref{eq:kinetic}) via a stochastic
model. Later the kinetic system
(\ref{eq:kinetic})-(\ref{eq:S_Poisson}) was formulated in
\cite{ODA88}. Othmer and Hillen in \cite{HO00, OH02} studied the
formal diffusion limit of this system by moment expansions. In
\cite{CMPS04} Chalub et al. proved that, in three dimensions with
suitable assumptions on the turning kernel $T_\eps$, the solution
to the kinetic system (\ref{eq:kinetic})-(\ref{eq:S_Poisson})
globally exists for any initial total mass. Moreover, they gave a
rigorous proof of the convergence to the macroscopic limit for all
time invervals of existence of the limiting Keller-Segel system.
In the following work \cite{HKS05, HKS05b} the results are
extended to more general cases.

The first task of this work is to design an Asymptotic-Preserving
(AP) method for the kinetic system
(\ref{eq:kinetic})-(\ref{eq:S_log}) in the one dimensional case.
The ``Asymptotic Preserving'' scheme, introduced by \cite{Jin}, is
a suitable numerical method for the kinetic equation in such a way
that letting $\eps\to0$ with the mesh size and time step fixed,
the scheme becomes a scheme for the limiting system. In this case,
this means as $\eps\to0$, the moments of the solution to the
kinetic system (\ref{eq:kinetic})-(\ref{eq:S_log}) solve the
modified Keller-Segel system
(\ref{eq:KellerSegel})-(\ref{eq:S_log}) automatically. We refer to
\cite{Jin_APreview} and the references therein for a detailed
review on AP scheme.

As mentioned above, the solution to the modified Keller-Segel
system (\ref{eq:KellerSegel})-(\ref{eq:S_log}) over the critical
mass blows up in finite time $t^*$. But the global solution to the
kinetic system (\ref{eq:kinetic})-(\ref{eq:S_log}) exists. An
attractive idea is that, the solution of the kinetic system after
$t^*$ might be a good candidate to describe the behavior of
solution to the Keller-Segel system after blow-up. There are other
mechanisms of regularizing the Patlak-Keller-Segel system after
blow-up by saturating the effect of the chemoattractant,
introducing volume effects or cross-diffusion in the macroscopic
system, see \cite{CHJ11} for a discussion. However, to regularize
it by the kinetic system is more appealing since the kinetic
system is directly derived from interacting particle systems
models.

Bournaveas and Calvez \cite{BC09} studied the kinetic system
(\ref{eq:kinetic})-(\ref{eq:S_Poisson}) with the local turning
kernel $T_\eps$ defined below in (\ref{def:BC_kernel}). They
showed that, in the two dimensional spherically symmetric
coordinate, the solution exists below a critical mass $m_c$ and
blows up over a critical mass $M_c$. Two questions remains
unsolved. First, the large-time behavior of the solutions for
subcritical cases is not clear yet. Second, the two thresholds do
not match, i.e. $m_c<M_c$. Moreover, the appearance of blow-up in
the one dimensional case has not been clarified yet, see
\cite{S11}. The second aim of this work is to investigate these
questions through numerical simulation.

The outline of the paper is as follows. In section
\ref{sec:kinetic} the kinetic models we will work on are
described. Section \ref{sec:macrolimit} gives the macroscopic
limit of each model. Section \ref{sec:scheme} describes the AP
schemes for these models. Finally, in section
\ref{sec:numericresult} we present the simulation experiments with
our schemes and draw some conclusions.


\section{The kinetic models}
\label{sec:kinetic}

In this section we briefly describe the kinetic models we used in
our simulations. The chemoattractant $S$ is always given by the
log kernel convolution (\ref{eq:S_log}).

The turning kernel $T_\eps$ in the kinetic equation
(\ref{eq:kinetic}) needs to be specified. The turning kernel,
$T_\eps = T_\eps[S](t,x,v,v')\ge0$, measures the probability of
velocity jump of cells from $v'$ to $v$. To derive the
Keller-Segel equation (\ref{eq:KellerSegel}) as $\eps\to0$, one
has to incorporate both $O(1)$ and $O(\eps)$ scale into $T_\eps$.
In the following work, we consider $T_\eps$ in the form
$$
T_\eps[S](t,x,v,v') = F(v) + \eps T_1 + O(\eps^2).
$$
Here $F(v)\ge0$ gives the equilibrium of velocity distribution
when there is no directional preference. It is natural to assume
\begin{equation}
  \label{def:Equi}
  \left\{
  \begin{aligned}
    \int_V F(v) \ud v = 1  \\
    F(v) = F(\abs{v})
  \end{aligned}
  \right. \, .
\end{equation}
$T_1$ characterizes the directional preference. We can assume
$T_1\ge0$ since we are considering the reaction of cells to
chemoattractant. The cells have a larger probability to jump to a
preferred direction.

\subsection{1D nonlocal Model}

\label{sec:kinetic_nonlocal}

Now we employ the nonlinear kernel introduced in \cite{CMPS04},
 \begin{equation*}
T_\eps=T_\eps[S](t,x,v,v')=\alpha_+\psi(S(t,x),S(t,x+\eps v)) +
\alpha_-\psi(S(t,x),S(t,x-\eps v')).
 \end{equation*}

The first term means the cell decides a new direction to move
based on the detection of current environment and probable new
location, while the second term gives the influence of past memory
on the choice of the new direction. For simplicity we neglect the
second term. We take $\alpha_+ = 1$, $\alpha_- = 0$, and we
consider
\begin{equation}
 \label{def:CJY_kernel}
 \psi(S,\tilde{S}) = F(v)+(\tilde{S}-S)_+.
\end{equation}
Again $F=F(|v|)$ is the equilibrium function satisfying
(\ref{def:Equi}). For simplicity we introduce the notation
\begin{equation*}
  \label{def:deltaS}
\delta^{\eps} S (x,v) = \left( S(x+\eps v) - S(x)\right)_+.
\end{equation*}
Note that $\delta^{\eps} S=O(\eps)$. Then the kinetic equation
(\ref{eq:kinetic}) reads
\begin{equation}
\label{eq:CJY1d}
 \eps \frac{{\partial f  }}{{\partial
t}} +  v  \frac{{\partial f  }}{{\partial
x}}   = \frac{1}{\eps}\left( \left(F(v) +  \delta^{\eps} S(x,v) \right)\rho - \left(1+\int \delta^{\eps} S(x,v') \ud v'\right)f\right).
\end{equation}
where $x\in\Omega=[-x_{\max},x_{\max}]$, $v\in V =
[-v_{\max},v_{\max}]$, $|V| = \mbox{Vol}(V)$. We impose the
initial conditions,
\begin{subequations}
\begin{align}
\label{eq:inival_f}
 f  \left( {0, x,v,} \right) & = f^I \left( {x,v} \right) \ge 0, \\
\label{eq:inival_S}
S  \left( {0,x} \right) & = S^I \left( x \right) \ge 0,
\end{align}
\end{subequations}
and reflection boundary condition for $f$, Neumann boundary
condition for $S$,
\begin{subequations}
\begin{align}
\label{eq:bdry_f}
f  \left( {t, \pm x_{\max},v} \right) &= f  \left( {t, \pm x_{\max},-v} \right), \\
\label{eq:bdry_S}
\partial_x S|_{x=\pm x_{\max}} &= 0.
\end{align}
\end{subequations}
Equations (\ref{eq:CJY1d})-(\ref{eq:S_log}) gives the nonlocal
model in one dimension. The global solution exists, regardless of
the total mass.

\subsection{1D local Model}

In this section we summarize a turning kernel introduced by
Bournaveas and Calvez \cite{BC09}. Let
\begin{equation}
  \label{def:BC_kernel}
  \ds T_\eps=T_\eps[S](t,x,v,v')=F(v) + \eps \left(v\cdot\nabla S(x)\right)_+,
\end{equation}
where $F=F(|v|)$ is the equilibrium function satisfying
(\ref{def:Equi}). Then the kinetic equation (\ref{eq:kinetic}) in
one dimension reads
\begin{equation}
\label{eq:BC1d}
 \ds \eps \frac{{\partial f  }}{{\partial
t}} +  v  \frac{{\partial f  }}{{\partial
x}}   = \frac{1}{\eps}\left( \left(F(v) +  \eps \left(v\cdot\nabla S\right)_+ \right)\rho - \left(1+ c_1\eps\abs{\nabla S} \right)f\right),
\end{equation}
with $\ds c_1 = \int_V \left(v\cdot\nabla S \ / \abs{\nabla
S}\right)_+ \ud v = \frac{1}{2} \int_V \abs{v} \ud v.$ The same
initial conditions (\ref{eq:inival_f}) (\ref{eq:inival_S}), and
boundary conditions (\ref{eq:bdry_f}), (\ref{eq:bdry_S}) are
applied.

\begin{remark}
If we expand the kernel (\ref{def:CJY_kernel}) $ T   = T_0  + \eps
T_1  + O\left( {\eps ^2 } \right)$, and drop the terms higher than
first order, we get exactly the kernel (\ref{def:BC_kernel}).
\end{remark}

\subsection{2D local Model}

Plug the local kernel (\ref{def:BC_kernel}) into the kinetic
equation (\ref{eq:kinetic}) in 2D, one obtains,
\begin{equation}
  \label{eq:BC2d}
  \left\{  \begin{array}{ll}
  \ds \eps \partial_t f + v \cdot \nabla_x f = \frac{1}{\eps} \left( \rho F(v) - f_{\eps} + \eps \rho (v\cdot \nabla S)_+ - \eps c_2 |\nabla S| f \right), & v\in V = B(0,v_{\max}), \,\,x\in \R^2 \\
  -\nabla^2 S = \rho
  \end{array} \right.
\end{equation}
with $c_2 = \int_V \left(v \cdot \nabla S / \abs{\nabla S}\right)_+ \ud v $. The
equilibrium $F(v)$ satisfies (\ref{def:Equi}).

Now we assume the initial data $f^I$ is spherically symmetric:
$f^I(\Theta x, \Theta v) = f^I (x,v)$, for any rotation $\Theta$.
Then the solution $f$ to (\ref{eq:BC2d}) remains spherically
symmetric. Denote
\begin{equation*}
  \ds r = |x| \in [0, \infty), \quad \omega = |v| \in [0,v_{\max} ), \quad \theta = \cos^{-1}\frac{x\cdot v}{r\omega} \in [0,\pi].
\end{equation*}
Then $f$ is a function of $r$, $\omega$, $\theta$ and $t$. Let
\begin{equation*}
  \ds h(t,r,\omega,\theta) = r f(t, x|_{(r,0)},v|_{(\omega\cos\theta, \omega\sin\theta)}) = r f(t,x|_{(r\cos\theta, r\sin\theta)},v|_{(\omega,0)}),
\end{equation*}
\begin{equation*}
  \ds \tilde{\rho}(t,r) = r\rho(t,r).
\end{equation*}
Then the density is given by
\begin{equation*}
  \ds \tilde{\rho}(t,r) = r\rho(t,r) = r\int_V f(t,|x|=r,v)\ud v = 2\iint_{0\le \omega \le v_{\max}, 0\le \theta \le \pi} \omega h \ud \theta \ud \omega,
\end{equation*}
and the total mass is
\begin{equation*}
  \ds M = 2\pi\int_0^{\infty} \tilde{\rho}(t,r) \ud r = 4\pi\iiint_{0\le \omega \le v_{\max}, 0\le \theta \le \pi,0\le r} \omega h \ud \theta \ud \omega \ud r.
\end{equation*}

Then equation (\ref{eq:BC2d}) can be rewritten as
\begin{equation}
\label{eq:BC2d_rad}
\left\{\begin{array}{l}
  \ds \partial_t h + \frac{\omega}{\eps }\left(\partial_r\left(\cos \theta h\right)-\partial_{\theta}\left(\frac{\sin \theta}{r} h\right)\right) = \frac{1}{\eps^2} ( \tilde{\rho} F(\omega) - h) + \frac{1}{\eps} \left(\omega \tilde{\rho} (\cos\theta\partial_r S)_+ - c_2 |\partial_r S| h\right),
  \\[3mm]
  -\partial_r \left(r \partial_r S\right) = \tilde{\rho},
\end{array} \right.
\end{equation}
where
$$
c_2 = 2\int_0^{\pi} \int_0^{v_{\max}} \omega^2 \left(\cos\theta \partial_r S / \abs{\partial_r S}\right)_+ \ud \omega \ud \theta = \int_0^{\pi} \int_0^{v_{\max}} \omega^2 \abs{\cos\theta} \ud \omega \ud \theta = \frac{2v_{\max}^3}{3},
$$
and $F(\omega)$ satisfies,
\begin{equation*}
  2\pi \int_0^{v_{\max}} \omega F(\omega) \ud \omega = 1.
\end{equation*}
From the second equation of (\ref{eq:BC2d_rad}), $\partial_r S$ can be computed by
\begin{equation*}
  \partial_r S = -\frac{1}{r}\int_0^r \tilde{\rho} \ud r.
\end{equation*}
We impose the boundary conditions,
\begin{equation*}
  \begin{aligned}
  h(0,\omega,\theta) &= h(0,\omega,\pi - \theta), \\
  \partial_r h(r_{\max},\omega,\theta) &= 0.
  \end{aligned}
\end{equation*}
It has been shown that the solution blows up in finite time in
\cite{BC09} for $M$ large enough while existing globally for small
enough mass.

\section{The macroscopic limits}
\label{sec:macrolimit}
\subsection{1D model: local and nonlocal}

The local kinetic model (\ref{eq:BC1d}) and the nonlocal kinetic model (\ref{eq:CJY1d}) give the same asymptotic limit as $\eps \to 0$.
We apply the Hilbert expansion into (\ref{eq:BC1d}) and (\ref{eq:CJY1d}), and match terms of the
same order in $\eps$. The classical Keller-Segel system can be derived for $\ds\rho(t,x) = \int_V f(t,x,v) \ud v$,
\begin{equation}
  \label{eq:KellerSegel_1d}
  \partial_t\rho = \partial_x\left(D\partial_x\rho - \chi\rho\partial_x S\right)
\end{equation}
where
\begin{equation}
  \label{eq:Dchi}
  D=\int_V|v|^2F(v)\ud v,
\qquad \mbox{ and } \qquad
  \chi=\frac{1}{2}\int_V|v|^2\ud v.
\end{equation}
We refer to \cite{CMPS04} for the details. Besides, by taking the
moments of the kinetic equation (\ref{eq:kinetic}), one has
\[
\ds\partial_t \rho + \nabla_x\cdot\left(\frac{1}{\eps} \int_V v f \ud v\right) = 0.
\]

Compare with (\ref{eq:KellerSegel_1d}), one has
\[
D\partial_x\rho - \chi\rho\partial_x S = -\frac{1}{\eps} \int_V v f \ud v.
\]
While the reflection boundary condition (\ref{eq:bdry_f}) leads to
\[
\int_V v f \ud v = 0, \quad \mbox{ at } x = \pm x_{\max}.
\]
One arrives at the general boundary condition (for mass preservation) for Keller-Segel model
\begin{equation*}
  D\partial_x\rho - \chi\rho\partial_x S = 0, \quad \mbox{ at } x=\pm x_{\max}.
\end{equation*}
Furthermore, the Neumann boundary condition for $\rho$ is derived under the condition (\ref{eq:bdry_S}).

\subsection{2D spherical symmetric local model}

For the reduced spherical symmetric kinetic system
(\ref{eq:BC2d_rad}), one can derive the similar asymptotic limit
for $\tilde \rho$ as $\eps \to 0$,
\begin{equation}
  \label{eq:KellerSegel_tilde}
  \partial_t \tilde\rho = \partial_r\left(Dr\partial_r\frac{\tilde\rho}{r} - \chi\tilde\rho\partial_r S\right)
\end{equation}
with
\begin{equation*}
  -\partial_r \left(r \partial_r S\right) = \tilde{\rho}
\end{equation*}
where $r\in[0,r_{\max}]$ with
\begin{equation*}
  D=2\int_0^{\pi} \int_0^{v_{\max}} \omega^3\cos^2\theta F(\omega)\ud \omega \ud \theta = \pi\int_0^{v_{\max}} \omega^3 F(\omega)\ud \omega
\end{equation*}
and
\begin{equation*}
  \chi=\int_0^{\pi} \int_0^{v_{\max}} \omega^3\cos^2\theta \ud \omega \ud \theta = \frac{\pi
  v_{\max}^4}{8}.
\end{equation*}
The Neumann boundary condition is derived,
\begin{equation*}
\begin{aligned}
  &\partial_r S = 0, &\quad\mbox{ at } r=0, r_{\max},\\
  &\partial_r\frac{\tilde\rho}{r} = 0, &\quad\mbox{ at } r=r_{\max}.
\end{aligned}
\end{equation*}


\section{The Numerical Scheme}
\label{sec:scheme}

In this section we present the numerical method for the kinetic
system ({\ref{eq:CJY1d}})-({\ref{eq:bdry_S}}) by the even and odd
parity formalism, which has been successfully applied to the
diffusive limit of linear transport equations in \cite{JPT01}.

\subsection{Odd and Even Parity}

We introduce the operator
\begin{equation*}
\begin{aligned}
  R[f](x,v) &= \frac{1}{2}\left(f(x,v)+f(x,-v)\right), \\
  J[f](x,v) &= \frac{1}{2\eps}\left(f(x,v)-f(x,-v)\right).
\end{aligned}
\end{equation*}
The new functions $R[f]$ and $J[f]$ are defined in $\Omega\times
V^+\times\mathbb{R}_+$, where $V^+ = \left\{v\in V | v\ge0
\right\}$. Then we take
\begin{eqnarray*}
 r(x,v) = R[f] = \frac{1}{2}\left(f(x,v)+f(x,-v)\right), & & j(x,v) = J[f] = \frac{1}{2\eps}\left(f(x,v)-f(x,-v)\right).
\end{eqnarray*}
We can recover $f$ from $r$ and $j$,
\begin{equation*}
 f(x,v) = \left\{
\begin{array}{l l}
 r(x,v) + \eps j(x,v),  &\quad \mbox{ if } v\ge0, \\
 r(x,-v) - \eps j(x,-v), &\quad \mbox{ if } v<0. \\
\end{array} \right.
\end{equation*}

\subsubsection{1D nonlocal model}

Now we describe the odd and even decomposition for the nonlocal model. Since $f(x,-v)$ satisfies
\begin{equation*}
 \eps \frac{{\partial f(x,-v)  }}{{\partial
t}} -  v  \frac{{\partial f(x,-v)  }}{{\partial
x}}   = \frac{1}{\eps}\left( \left(F(v) +  \delta^{\eps} S(x,-v) \right)\rho - \left(1+\int \delta^{\eps} S(x,v') \ud v'\right)f(x,-v)\right),
\end{equation*}
we obtain,
\begin{eqnarray*}
  \partial_t r + v \partial_x j &=&
  \frac{1}{\eps^2}\left(\left(F(v)+R[\delta^{\eps} S]\right)\rho - \left(1+<\delta^{\eps} S> \right)r\right), \\
  \partial_t j + \frac{1}{\eps^2}v \partial_x r &=&
  \frac{1}{\eps^2}\left(J[\delta^{\eps} S]\rho - \left(1+ <\delta^{\eps} S> \right)j\right),
\end{eqnarray*}
where $<\delta^{\eps} S> = \int \delta^{\eps} S(x, v') \ud v'$, and
\begin{equation}
\label{fmla:rho}
\rho=\int_V f \ud v =
\int_{V^+}\left(f(x,v)+f(x,-v))\right)\ud v=2\int_{V^+}r(x,v)\ud v.
\end{equation}
We assume $\eps \le 1$, then we can rewrite the equations for $r$ and $j$ as,
\begin{equation}
\label{eq:rj}
  \begin{aligned}
  \partial_t r + v \partial_x j &=
  \frac{1}{\eps^2}\left(\left(F(v)+R[\delta^{\eps} S]\right)\rho - \left(1+<\delta^{\eps} S> \right)r\right), \\
  \partial_t j + v \partial_x r &=
  \frac{1}{\eps^2}\left(J[\delta^{\eps} S]\rho - \left(1+<\delta^{\eps} S> \right)j\right) + \left(1-\frac{1}{\eps^2}\right)v \partial_x r.
  \end{aligned}
\end{equation}

Let us finish with the boundary conditions for new variables $r$
and $j$. From the reflection boundary condition (\ref{eq:bdry_f})
and the definition of $r$ and $j$, we can easily get the boundary
condition for $j$,
\begin{equation}
\label{eq:bdry_j} j(t,x_b,v) = \frac{1}{2\eps}\left ( f(t,x_b,v) -
f(t,x_b,-v)\right) = 0.
\end{equation}
where $x_b=\pm x_{\max}$. Then plug this into the second equation
of (\ref{eq:rj}), we derive ,
\[
v\partial_x r|_{x=x_b}= J[\delta^{\eps} S]\rho |_{x=x_b} = \frac{1}{2\eps}(\delta^{\eps} S(x_b, v)-\delta^{\eps} S(x_b,-v)) \rho.
\]
Now from the Neumann boundary condition (\ref{eq:bdry_S}), we have
$S(x_b+\eps v)=S(x_b-\eps v)$, i.e. $\delta^{\eps} S(x_b,
v)=\delta^{\eps} S(x_b,-v)$. So the boundary condition for $r$ is
\begin{equation}
\label{eq:bdry_r}
\partial_x r|_{x=x_b}= 0.
\end{equation}

Now, we can propose an asymptotic preserving method for the one
dimensional nonlocal model. The idea can be applied to the other
models straightforwardly.

As $\eps$ approaches $0$, we can derive the Keller-Segel
equation (\ref{eq:KellerSegel}) asymptotically from the kinetic
equation (\ref{eq:CJY1d}). So a natural requirement is the numerical
method for (\ref{eq:CJY1d}) should discretize its macroscopic limit
as $\eps \to 0$. We give an AP method following \cite{JP00}.

We can employ an operator splitting method on
(\ref{eq:rj}). First the (stiff) source part is
solved by the implicit Euler method,
\begin{equation}
\label{num:source}
\begin{aligned}
  \partial_t r &=
  \frac{1}{\eps^2}\left(\left(F(v)+R[\delta^{\eps} S]\right)\rho - \left(1+<\delta^{\eps} S> \right)r\right), \\
  \partial_t j &=
  \frac{1}{\eps^2}\left(J[\delta^{\eps} S]\rho - \left(1+<\delta^{\eps} S> \right)j\right) + \left(1-\frac{1}{\eps^2}\right)v \partial_x r.
\end{aligned}
\end{equation}
Then the transport part can be solved by an explicit method (e.g. upwind scheme),
\begin{equation}
  \label{num:transport}
  \begin{aligned}
    \partial_t r + v \partial_x j &= 0, \\
  \partial_t j + v \partial_x r &= 0.
  \end{aligned}
\end{equation}
We can check that the method described above is AP easily. As
$\eps\to0$ the leading term in $\eps$ of
(\ref{num:source}) gives,
\begin{eqnarray*}
r &=&\frac{F(v)+R[\delta^{\eps} S]}{1+<\delta^{\eps} S>}\rho= \rho F(v)+O(\eps), \\
j &=&\frac{J[\delta^{\eps} S]\rho-v\partial_x r}{1+<\delta^{\eps}
S>}= v\left(\frac{1}{2}\partial_x S\rho-\partial_x r\right) +
O(\eps).
\end{eqnarray*}
Plug into the first equation of (\ref{num:transport}) and
integrate over $V^+$, we get exactly the Keller-Segel equation
(\ref{eq:KellerSegel}) with $D$ and $\chi$ given by
(\ref{eq:Dchi}).


\subsubsection{1D local model}
For the one dimensional local model, we obtain,
\begin{eqnarray*}
  \partial_t r + v \partial_x j &=&
  \frac{1}{\eps^2}\left(\left(F(v)+ \frac{\eps}{2}\abs{v\partial_x S}\right)\rho - (1+c_1 \eps \abs{\partial_x S} )r\right), \\
  \partial_t j + \frac{1}{\eps^2}v \partial_x r &=&
  \frac{1}{\eps^2}\left(\frac{1}{2}v\partial_x S \rho - (1+c_1 \eps \abs{\partial_x S} )j\right).
\end{eqnarray*}
The remaining work is similar to the nonlocal model. The boundary conditions for for new variables $r$ and $j$ are also given by (\ref{eq:bdry_r}) and (\ref{eq:bdry_j}). The AP scheme can be designed similarly.

\subsubsection{2D spherical symmetric local model}

Denote
\begin{equation*}
  R = \frac{1}{2} (h(r,\omega,\theta) + h(r,\omega,\pi - \theta)), \hskip 1cm J = \frac{1}{2\eps} (h(r,\omega,\theta) - h(r,\omega,\pi - \theta)).
\end{equation*}
$R$ and $J$ are functions defined on
$[0,r_{\max}]\times[0,v_{\max}]\times[0,\pi]$. Then we can write
the equations for $R$ and $J$,
\begin{equation*}
\begin{aligned}
  \partial_t R + \omega \left(\partial_r (\cos\theta J) - \partial_{\theta} \left(\frac{\sin\theta}{r}J\right)\right) &= \frac{1}{\eps^2} \left(\tilde{\rho} F - R \right) + \frac{1}{\eps}\left( \frac{1}{2} \omega\tilde{\rho}\cos\theta |\partial_r S| - c_2 |\partial_r S| R \right), \\
  \partial_t J + \frac{\omega}{\eps^2} \left(\partial_r (\cos\theta R) - \partial_{\theta} \left(\frac{\sin\theta}{r}R\right)\right) &= -\frac{1}{\eps^2} J  + \frac{1}{\eps}\left( \frac{1}{2\eps} \omega\tilde{\rho}\cos\theta \partial_r S - c_2 |\partial_r S| J \right).
\end{aligned}
\end{equation*}
Now we split this into two steps. First the collision step can be solved implicitly,
\begin{equation*}
\begin{aligned}
  \partial_t R &= \frac{1}{\eps^2} \left(\tilde{\rho} F - R \right) + \frac{1}{\eps}\left( \frac{1}{2} \omega\tilde{\rho}\cos\theta |\partial_r S| - c_2 |\partial_r S| R \right), \\
  \partial_t J &= -\frac{1}{\eps^2} J  + \frac{1}{\eps}\left( \frac{1}{2\eps} \omega\tilde{\rho}\cos\theta \partial_r S - c_2 |\partial_r S| J \right) + \left(1-\frac{1}{\eps^2}\right)\omega \left(\partial_r (\cos\theta R) - \partial_{\theta} \left(\frac{\sin\theta}{r}R\right)\right).
\end{aligned}
\end{equation*}
Then the transport part can be solved explicitly,
\begin{equation}
\label{eq:BC2d_rad_RJtransport}
\begin{aligned}
  \partial_t R + \omega \left(\partial_r (\cos\theta J) - \partial_{\theta} \left(\frac{\sin\theta}{r}J\right)\right) = 0, \\
  \partial_t J + \omega \left(\partial_r (\cos\theta R) - \partial_{\theta} \left(\frac{\sin\theta}{r}R\right)\right) = 0.
\end{aligned}
\end{equation}
One can show $\tilde{\rho}$ solves the macroscopic limit (\ref{eq:KellerSegel_tilde}) as $\eps \to 0$.

\subsection{Time and space discretization}

\subsubsection{1D system: nonlocal and local}

In this section we give full discretization for the 1D nonlocal
kinetic system. The discretization for the 1D local system is
similar. First we solve (\ref{num:source}) by implicit Euler
method,
\begin{equation*}
  \begin{aligned}
    \frac{r^* -r^n}{\Delta t} &=\frac{1}{\eps^2}\left( (F(v)+R[\delta^{\eps} S^*])\rho^* - (1+<\delta^{\eps} S^*>)r^*\right), \\
  \frac{j^* -j^n}{\Delta t} &=\frac{1}{\eps^2}\left(J[\delta^{\eps} S^*]\rho^*-(1+<\delta^{\eps} S^*>)j^* \right) + \left(1-\frac{1}{\eps^2}  \right) v \partial_x^{(c)} r^*.
  \end{aligned}
\end{equation*}
where $r^n$ and $j^n$ are the numerical solutions of $r$ and $j$
at time $t_n$. $\partial_x^{(c)}$ is the central difference
discretization of $\partial_x$. And a linear interpolation is used
to get $S(x+\eps v)$ and $S(x-\eps v)$. $<\delta^{\eps} S>$ is
evaluated by applying the trapezoidal rule on the interpolated
value.

If we integrate (\ref{num:source}) over $V^+$, we can get $\partial_t \rho =0$. So
\begin{equation*}
  \rho^* = \rho^n, \hskip 1cm S^* = S^n.
\end{equation*}
Therefore $r^*$ and $j^*$ can be solved \textit{explicitly},
\begin{equation*}
  \begin{aligned}
  r^* &= \frac{\eps^2 r^n + \Delta t (F(v)+ R[\delta^{\eps} S^n])\rho^n}{\eps^2 + \Delta t (1+<\delta^{\eps} S^n>)},\\
  j^* &= \frac{\eps^2 j^n + \Delta t \left(J[\delta^{\eps} S^n]\rho^n + (\eps^2-1)v\partial_x^{(c)} r^*\right)}{\eps^2 + \Delta t (1+<\delta^{\eps} S^n> )}.
  \end{aligned}
\end{equation*}
Then we apply a first-order upwind scheme  on
(\ref{num:transport}) to get $r^{n+1}$, $j^{n+1}$,
\begin{equation*}
  \begin{aligned}
 \frac{r_i^{n+1}-r_i^n}{\Delta t} + \frac{v}{2\Delta x}\left(j_{i+1}-j_{i-1}\right) &= \frac{v}{2\Delta x}\left(r_{i+1}-2r_i+r_{i-1}\right), \\
 \frac{j_i^{n+1}-j_i^n}{\Delta t} + \frac{v}{2\Delta x}\left(r_{i+1}-r_{i-1}\right) &= \frac{v}{2\Delta x}\left(j_{i+1}-2j_i+j_{i-1}\right).
  \end{aligned}
\end{equation*}

As $\eps\to0$, a simple computation shows that this scheme leads
to (after integration over $V^+$),
\begin{equation*}
   \frac{\rho^{n+1}-\rho^n}{\Delta t}= \partial_x^{(c)}\left(D\partial_x^{(c)}\rho^n -
   \chi\rho^n\partial_x^{(c)}S^n\right) + \Delta x C(V)\partial_{xx}^{(c)}\rho^n.
\end{equation*}
with $D$ and $\chi$ given by (\ref{eq:Dchi}), $\rho$ given by
(\ref{fmla:rho}). Here $\partial_{xx}^{(c)}$ is the general
three-point central difference of $\partial_{xx}$. $C(V)$ is a
constant which only depends on the velocity space $V$. In the case
of $V=[-v_{\max},v_{\max}]$, $C(V)=v_{\max}/4$.

\subsubsection{2D spherical symmetric local model}

We describe a first order discretization for the transport part
(\ref{eq:BC2d_rad_RJtransport}) of spherical symmetric system. By
introducing
\begin{equation*}
  P(x,v) = \frac{1}{2} (R(x,v) + J(x,v)), \hskip 1cm   Q(x,v) = \frac{1}{2} (R(x,v) - J(x,v)),
\end{equation*}
one obtains
\begin{equation*}
\begin{aligned}
  \partial_t P + \omega \left(\partial_r (\cos\theta P) - \partial_{\theta} \left(\frac{\sin\theta}{r}P\right)\right) = 0, \\
  \partial_t Q - \omega \left(\partial_r (\cos\theta Q) - \partial_{\theta} \left(\frac{\sin\theta}{r}Q\right)\right) = 0.
\end{aligned}
\end{equation*}
We take the grid points at
\begin{equation*}
\begin{aligned}
    r_i &= (i-\frac{1}{2})\Delta r, & \hskip .3cm i = 1,\dots,N, \\
    \theta_j &= (j-\frac{1}{2})\Delta \theta, & \hskip .3cm j = 1,\dots,M.
\end{aligned}
\end{equation*}

 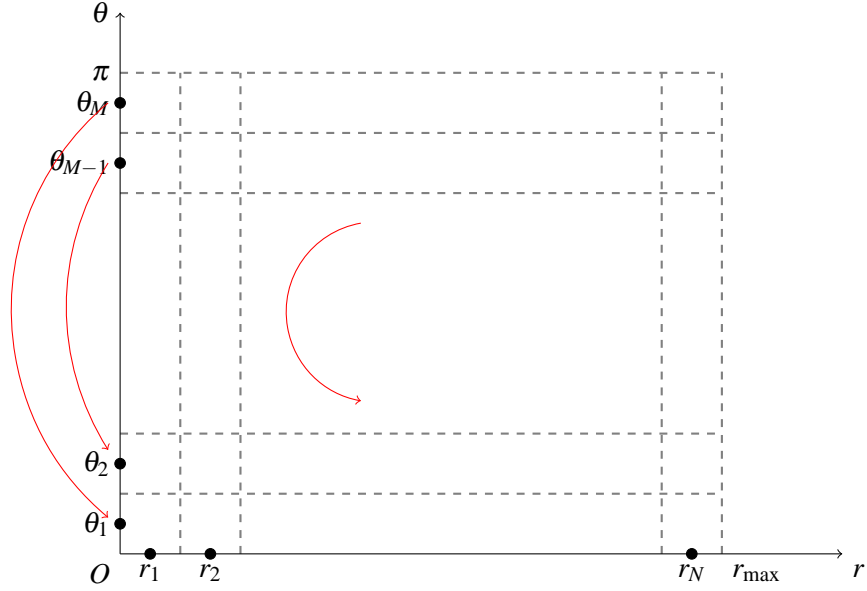
\begin{figure}
 \centering
\begin{tikzpicture}[scale=0.8]
\draw[thick,color=gray,dashed]  (0,1) -- (10,1);
\draw[thick,color=gray,dashed]  (0,2) -- (10,2);
\draw[thick,color=gray,dashed]  (0,6) -- (10,6);
\draw[thick,color=gray,dashed]  (0,7) -- (10,7);
\draw[thick,color=gray,dashed]  (0,8) -- (10,8);
\draw[thick,color=gray,dashed]  (1,0) -- (1,8);
\draw[thick,color=gray,dashed]  (2,0) -- (2,8);
\draw[thick,color=gray,dashed]  (9,0) -- (9,8);
\draw[thick,color=gray,dashed]  (10,0) -- (10,8);
\draw[->] (0,0) -- (12,0)
node[below right] {$r$};
\draw[->] (0,0) -- (0,9)
node[left] {$\theta$};
\draw[fill](.5,0)circle(2.5pt) node[below] {$r_1$};
\draw[fill](1.5,0)circle(2.5pt) node[below] {$r_2$};
\draw[fill](9.5,0)circle(2.5pt) node[below] {$r_N$};
\draw[fill](0,0.5)circle(2.5pt) node[left] {$\theta_1$};
\draw[fill](0,1.5)circle(2.5pt) node[left] {$\theta_2$};
\draw[fill](0,6.5)circle(2.5pt) node[left] {$\theta_{M-1}$};
\draw[fill](0,7.5)circle(2.5pt) node[left] {$\theta_M$};
\draw(0,0) node[below left] {$O$};
\draw(0,8) node[left] {$\pi$};
\draw(10,0) node[below right] {$r_{\max}$};
\draw[->,color=red](-0.2,7.5) arc (130:230:4.5cm);
\draw[->,color=red](-0.2,6.5) arc (148:212:4.5cm);
\draw[->,color=red](4,5.5) arc (100:260:1.5cm);
\end{tikzpicture}
\caption{Illustration of grid points and characteristic lines on $r-\theta$ plane.}
\label{fig:transportdemo}
\end{figure}

The grid points and the characteristic lines are shown in Figure
\ref{fig:transportdemo}. We can define the flux at each interface
according to the ``upwind'' value along the characteristic
direction. An analogous idea was successfully used in the case of
numerical simulation of the Boltzmann-Poisson system for
semiconductors in \cite{CGMS}. Here we give the detailed
discretization for $P$. The discretization for $Q$ is similar. For
each $\omega$,
\begin{equation*}
  \frac{P^{n+1}_{i,j} - P^{n}_{i,j}}{\Delta t} + \omega \left(\frac{F^n_{i+1/2,j} - F^n_{i-1/2,j}}{\Delta r}  + \frac{G^n_{i,j+1/2} - G^n_{i,j-1/2}}{\Delta \theta} \right) = 0,
\end{equation*}
where
\begin{equation*}
  F^n_{i+1/2,j} =
  \left\{
  \begin{aligned}
    &(\cos \theta _j) P^n_{i,j}, &\mbox{ if } j \le M/2,  \\
    &(\cos \theta _j) P^n_{i+1,j}, &\mbox{ if } j \ge M/2+1,
  \end{aligned}
  \right.
\end{equation*}
\begin{equation*}
  G^n_{i,j+1/2} = -\frac{\sin(\theta_j+\Delta \theta/2)}{r_i} P^n_{i,j} = -\frac{\sin(j \Delta \theta)}{r_i} P^n_{i,j},
\end{equation*}
with boundary
\begin{equation*}
  F^n_{1/2,j} = F^n_{1/2,M-j} = (\cos \theta _{M-j}) P^n_{1,M-j}, \mbox{ if } j \le M/2.
\end{equation*}

\subsection{Adaptive grids for solutions with blowup}
\label{sec:adaptivegrid} Now let us consider the case that the
initial total mass is large. It has been shown that the solutions
to the 2D spherical symmetric kinetic model (\ref{eq:BC2d_rad})
blow up in finite time. While for the one dimensional kinetic
local model (\ref{eq:BC1d}), although the strict analysis is
lacking \cite{S11}, the numerical results strongly suggest that
the solutions also blow up. In these cases, the convergence of the
schemes described above would be questionable when simulations are
performed on the fixed grids. See section for the verification of
convergence order.

Noting that in the blowup case, the magnitude of $\nabla_x S$
grows up dramatically as time evolves, which actually
characterizes another scale for the kinetic equation. To capture
this scale, one has to take $\Delta x\sim \frac{1}{\max |\nabla_x
S|}$. Hence the adaptive grids are needed. In numerical
simulation, we double the grid points once $||\nabla_x
S||_{\infty}$ is doubled. More exactly, we apply the following
algorithm.

\begin{algori}
\emph{(Adaptive grids)} \
\begin{algorithmic}
 \STATE $s = ||\nabla_x S_0||_{\infty}$
  \IF {$||\nabla_x S_n||_{\infty}\geq 2s$}
        \STATE $s \gets ||\nabla_x S_n||_{\infty}$.
        \STATE Double the grids, half the time step size $\Delta t$.
        \STATE Derive $f^n$ on the finer grids by interpolation and continue the time evolution.
  \ENDIF
\end{algorithmic}
\end{algori}

This simple idea works well when checking the blowup property and
determining the blowup time. A nonuniform refinement might be more
efficient in this problem since the mass is concentrated at some
separated points. However this is beyond the scope of this work
and is left for future study.

\subsection{A second order scheme for 1D model}
\label{sec:2ndorder}

One can derive a second order scheme without much difficulty for
the 1D (local or nonlocal) kinetic model. First the Strang
splitting is applied. One solves the stiff source part
(\ref{num:source}) over a time step $\frac{\Delta t}{2}$, then
solves the transport part (\ref{num:transport}) over $\Delta t$,
then solves the stiff source part (\ref{num:source}) over another
$\frac{\Delta t}{2}$.

Here we describe the second order scheme for the nonlocal model in
detail. It can be applied to the local model without any
difficulty. The stiff part (\ref{num:source}) can be solved
exactly in time. Noting that the cell density $\rho$ and the
chemical concentration $S$ are not changed in this step. Let
$$
\lambda = e^{-\frac{\Delta t (1+<\delta^{\eps}S^n>)}{2\eps^2}}\,,
$$
then after a time step $\Delta t/2$,
\begin{equation*}
  r^* =  \lambda r^n + \left(1-\lambda\right) \frac{F(v)+R[\delta^{\eps}S^n]}{1+<\delta^{\eps}S^n>}\rho^n,
\end{equation*}
\begin{equation*}
  j^* = \lambda j^n + \left(1-\lambda\right) \frac{J[\delta^{\eps}S^n]}{1+<\delta^{\eps}S^n>}\rho^n + v \partial_x^{(c)} \tilde r,
\end{equation*}
with
\begin{equation*}
  \tilde r = \left(1-\frac{1}{\eps^2}\right) \left(\frac{\Delta t}{2} \lambda r^n + \frac{F(v)+R[\delta^{\eps}S^n]}{1+<\delta^{\eps}S^n>}\rho^n \left( \frac{\eps^2  \left(1-\lambda\right) }{(1+<\delta^{\eps}S^n>)} - \frac{\Delta t}{2} \lambda
  \right)
  \right).
\end{equation*}
Here the computation of $\delta^{\eps}S^n$ is described in Section \ref{sec:deltaS}. $\partial_x^{(c)}$ is a central difference in $x$
direction.

Next a second-order TVD scheme is needed to solve the transport
part over a time step $\Delta t$,
\begin{equation*}
 \frac{r_i^{**}-r_i^*}{\Delta t} + \frac{v}{2\Delta x}\left(j_{i+1}^*-j_{i-1}^*\right) = \frac{v}{2\Delta x}\left(r_{i+1}^*-2r_i^*+r_{i-1}^*\right) -
\frac{v}{4}\left(1-\frac{v\Delta t }{\Delta x}\right)\left(\sigma^+_i-\sigma^+_{i-1}+\sigma^-_{i+1}-\sigma^-_i\right),
\end{equation*}
\begin{equation*}
 \frac{j_i^{**1}-j_i^*}{\Delta t} + \frac{v}{2\Delta x}\left(r_{i+1}^*-r_{i-1}^*\right) = \frac{v}{2\Delta x}\left(j_{i+1}^*-2j_i^*+j_{i-1}^*\right) -
\frac{v}{4}\left(1-\frac{v\Delta t }{\Delta x}\right)\left(\sigma^+_i-\sigma^+_{i-1}-\sigma^-_{i+1}+\sigma^-_i\right),
\end{equation*}
where
\begin{equation*}
\sigma^{\pm}_i=\frac{1}{\Delta x}\mbox{minmod}\left(r_{i+1}^* \pm j_{i+1}^*-r_i^* \mp j_i^*, \hskip .5cm r_i^* \pm j_i^*- r_{i-1}^* \mp j_{i-1}^*\right).
\end{equation*}
Finally we can update the density $\rho^{**}$ from $r^{**}$ and
$S^{**}$ is obtained. Then the stiff part is solved over another
$\Delta t/2$.

\subsection{The convolution}

Finally, we apply the FFT algorithm on the computation of
convolution (\ref{eq:S_log}) to get $S$. The singularity of
$\log|x|$ at $x=0$ makes the direct numerical integral difficult,
see \cite{CPS07}. Noting that $\log|x|$ belongs to $L^1$, we can
avoid this problem by taking Fourier transform first. In the
numerical simulation, we will always make $f$ to be compactly
supported in the computational region. It is not a problem about
the periodicity of boundary condition by extending the solution to
zero in a larger interval and by computing the Fourier transform
there. In this way we avoid any kind of aliasing.

\subsection{The computation of $\delta^{\eps}S$}
\label{sec:deltaS}

To obtain a higher accuracy in computing $\delta^{\eps}S$ in (\ref{def:deltaS}), a high order interpolation is needed to compute $S(x+\eps v)$. Here we apply the FFT based interpolation,
\[  \ds S(x_j+\eps v) = \frac{1}{N}\sum_k \hat{S}_k e^{-ik(x_j+\eps v)} = \frac{1}{N}\sum_k \left(\hat{S}_k e^{-ik \eps v}\right) e^{-ikx_j}, \]
where $\hat{S}$ is the discrete Fourier transform of $S$ on grids ${x_j}$.


\section{Numerical results}
\label{sec:numericresult}

\subsection{1D Nonlocal Model}

The first simulations are devoted to the one dimensional nonlocal model described in section \ref{sec:kinetic_nonlocal}.

The following simulations are set on $x\in\Omega=[-1,1]$, $v\in V = [-1, 1]$. We
take $N_v = 64$, which can provide good enough accuracy for
numerical simulations. The constant function in $V$ is chosen as the equilibrium
$$
\ds F(v) = \frac{1}{|V|} \mathbf{1}_{V}.
$$
In this setting, the critical mass of blow up for the limiting Keller-Segel system is
$$
M_{c} = 2\pi.
$$

The initial conditions in the
simulation are always given by,
\begin{equation}
\label{num:init_general}
 \rho^I(x) = C e^{-80x^2}, \qquad f^I_{\eps}(x,v) = \rho^I(x) F(v),
\end{equation}
where $C=C(M)$ is a constant determined by the total mass $M$.

As predicted in \cite{CMPS04}, in the supercritical case $M>M_c$,
the solutions to the kinetic system converge to that of the
Keller-Segel system only in finite time (before blow up time).
After that, the asymptotic limit is not valid anymore. To capture
the behavior of solutions to the kinetic system after that time,
one has to resolve the small scale $\eps$. Therefore in the
simulation, we need $\Delta x = O(\eps)$. While in the subcritical
case, as will be shown in the following sections, the asymptotic
limit seems to be valid over any time period. One can take $\Delta
x$ independent of $\eps$, as in a typical AP scheme \cite{Jin_APreview}.

As for the time step length $\Delta t$, the simulation results
suggest that, for the sake of stability, one needs $\Delta t =
\eps \Delta x/(2v_{\max})$ for long time simulation in the
supercritical case. While in the subcritical case $M<M_c$, as
$\eps\to0$, the diffusive nature of the Keller-Segel system
requires $\Delta t = \Delta x^2/2$. A general choice of $\Delta t$
would be
\[\ds \Delta t = \max\left\{\frac{\eps \Delta
x}{2v_{\max}}, \frac{\Delta x^2}{2}\right\}.\]

\subsubsection{Convergence order of numerical scheme}

In this section we test the convergence order of numerical scheme
described in Section \ref{sec:2ndorder}. We check the following error,
\begin{equation}
\label{num:err}
\ds  e_{\Delta x}(f) =  \frac{|| f_{\Delta x}(t) - f_{2\Delta x}(t) ||_1}{||f_{2\Delta x}(0)||_1}.
\end{equation}
This can be considered as an estimation of the relative error in
$l^1$ norm, where $f_h$ is the numerical solution computed from a
grid of size $\Delta x=\frac{x_{\max}-x_{\min}}{N_x}$. The
numerical scheme is said to be $k$-th order if $e_{\Delta x}\le C
{\Delta x}^k$. The computations are performed with $N_x = 100, 200, 400, 800, 1600, 1600$, $\Delta t = \frac{\Delta x^2}{2}$.

\begin{figure}
 \centering
 \subfigure[$M >M_c$.]  {\includegraphics[width=0.45\textwidth]{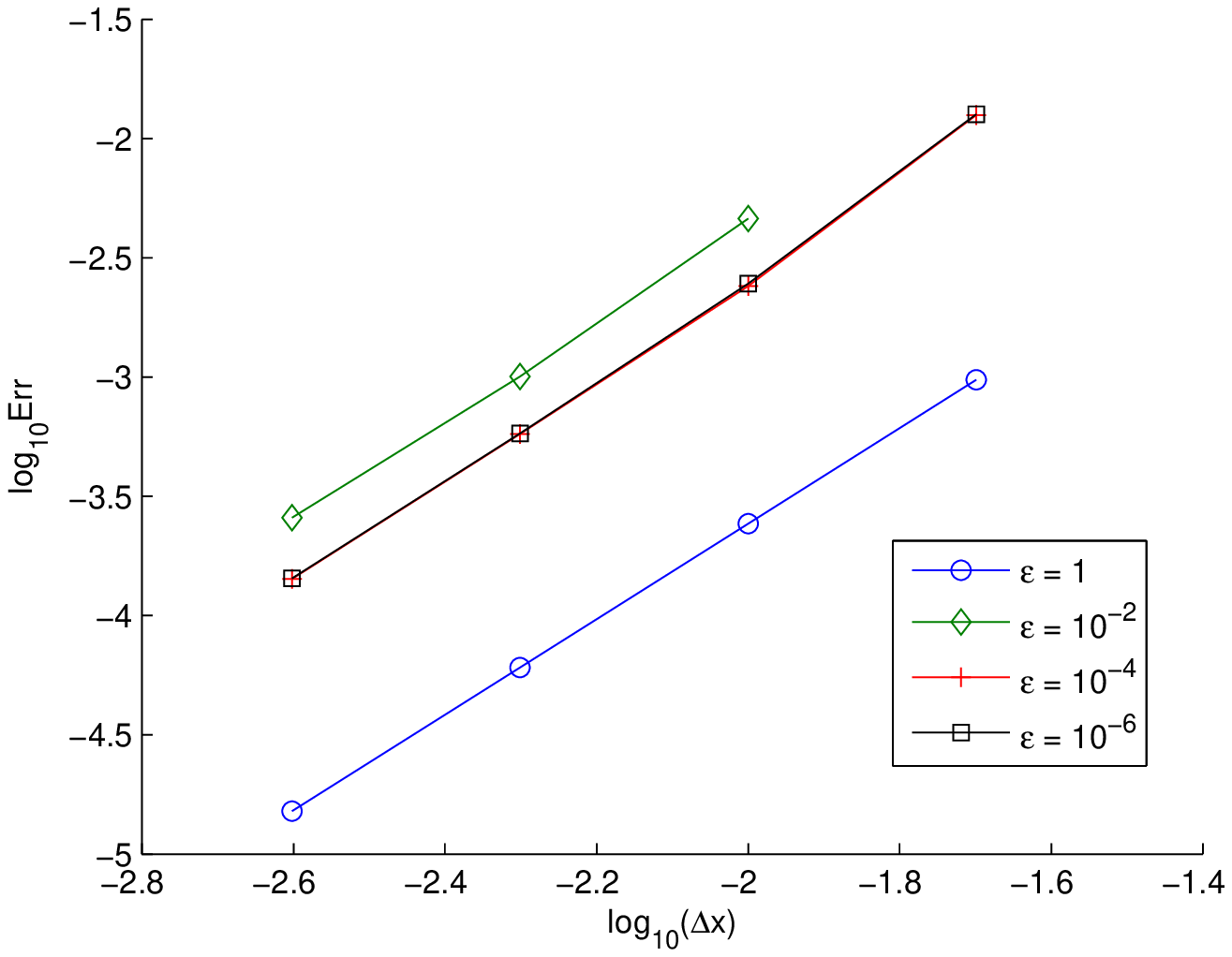}}
 \subfigure[$M <M_c$.]  {\includegraphics[width=0.45\textwidth]{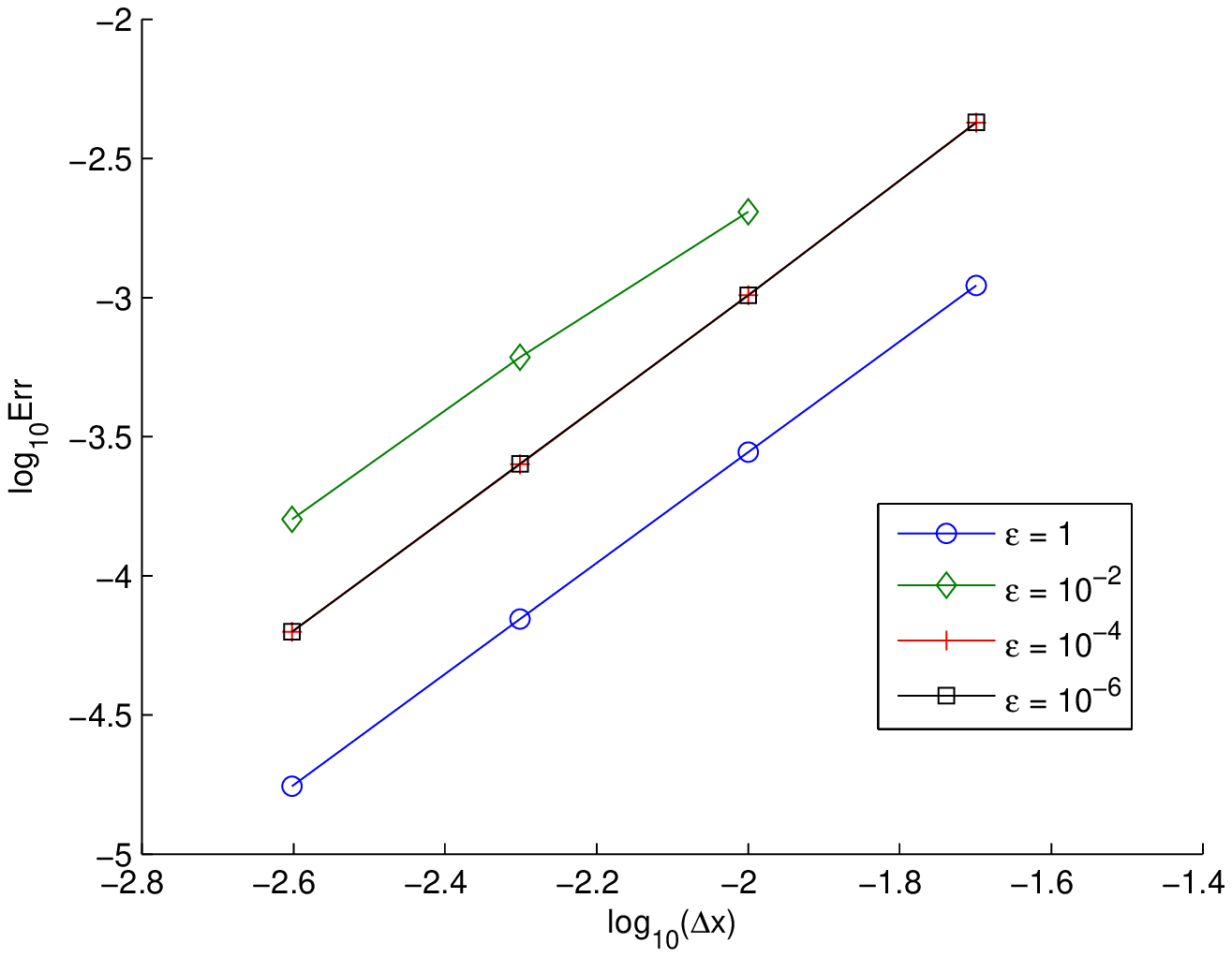}}
 \caption{{1D nonlocal model. The convergence order of scheme described in Section \ref{sec:2ndorder}, for different $\eps$. Left: $M=4\pi>M_c$, $t = 0.0025<t_b$; Right: $M=\pi<M_c$, $t=0.025$. }}
\label{fig:2ndorder}
\end{figure}

Figure \ref{fig:2ndorder} (a) shows the convergence results at $t = 0.0025$ for different $\eps$, with initial data given by (\ref{num:init_general}) and total mass $M=4\pi>M_c$. The solution to the Keller-Segel system with total mass $4\pi$ blows up at $t_b\approx0.0039$ (see next section).  The scheme shows  second order convergence (in $l^1$ norm) for supercritical mass before blow up time.

Figure \ref{fig:2ndorder} (b) shows the convergence results at $t = 0.025$ for different $\eps$, with initial data given by (\ref{num:init_general}) and subcritical mass $M=\pi<M_c$. The second order convergence (in $l^1$ norm) is observed for all $\eps$.

In conclusion, our scheme has second order convergence, uniformly in $\eps$. This is a common result for AP scheme, see \cite{GJL99}.

The above simulations are performed with the transport equation (\ref{num:transport}) solved by a Lax-Wendroff method. The use of the second order TVD method described in Section \ref{sec:2ndorder} shows a lower, but still uniform convergence in $\eps$.

\subsubsection{Global existence and finite time blow up}

Following the proof in \cite{CMPS04}, one can show that the
solution to the kinetic system (\ref{eq:CJY1d})-(\ref{eq:bdry_S})
is bounded on $[0,T]$, for any time $T$. However, the
Patlak-Keller-Segel system (\ref{eq:KellerSegel}) can present a
blow-up phenomenon in finite time, see \cite{CPS07}.

Now we take the initial data (\ref{num:init_general}) with supercritical mass $M = 4\pi >M_c = 2\pi$.
The blowup time $t_b\approx 0.0039$ from the numerical simulation.

\begin{figure}
 \centering
 \includegraphics[width=\textwidth]{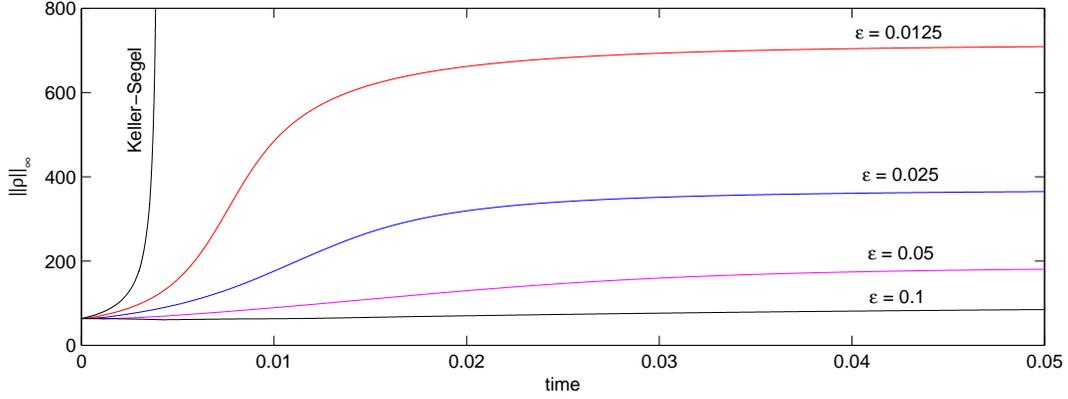}
 \caption{{1D nonlocal model. Time evolution of kinetic system and Keller-Segel model for supercritical mass $M=4\pi$.}}
\label{fig:kin_KS}
\end{figure}

Figure \ref{fig:kin_KS} shows the global boundedness of kinetic
system for different $\eps$ and the finite time blowup for the
corresponding Keller-Segel model, by drawing the time evolution of
the maximum value of $\rho_{\eps}$ and $\rho$.

\subsubsection{AP property: Convergence in $\eps$ at a finite time interval}

As mentioned before, it has been shown in \cite{CMPS04} that the solution of the kinetic system
can converge to that of the Patlak-Keller-Segel system weakly in a
finite time interval $[0, t^*]$ with $t^*$ small enough. Here we
numerically check this convergence. We use the same grid size $\Delta x = \frac{1}{2000}$ for
different $\eps$ in this subsection.

In this section, we use the notations $\rho_\eps$ and $f_\eps$ to for the solutions to the kinetic system, while $\rho_0$ for the Keller-Segel model.

\begin{figure}
\centering
\subfigure[$\rho_\eps(x)$.]{ \includegraphics[width=0.45\textwidth]{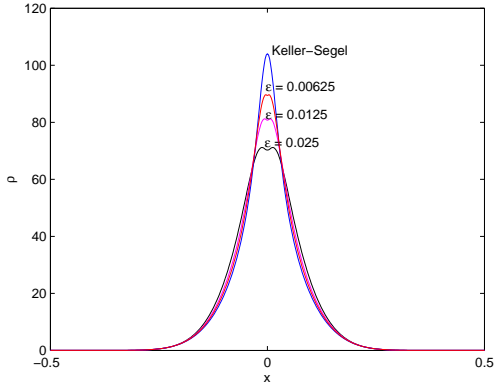}}
\subfigure[Convergence in $\eps$.]{ \includegraphics[width=0.45\textwidth]{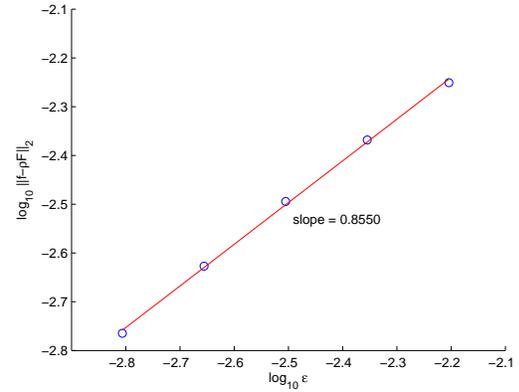}}
 \caption{{1D nonlocal model. The solutions of kinetic system and Keller-Segel system at time $t = 0.002$, before blow up time. The super-critical mass $M=4\pi$ is used. The left figure shows $\rho_{\eps}$ (for different $\eps$) and $\rho$. The right figure gives the convergence order of $||f_\eps(x,v) - \rho_\eps(x) F(v)||_2$.}}
\label{fig:conv_beforeBU}
\end{figure}

The supercritical case is studied first. We take initial data (\ref{num:init_general}) with $M=4\pi$.
Figure \ref{fig:conv_beforeBU} shows the
convergence of $f_{\eps}\to \rho_0 F$ as $\eps \to 0$ at time
$t=0.002<t_b$, where $t_b\approx0.0039$ is the blow up time of the limiting Keller-Segel model.
Figure \ref{fig:conv_beforeBU}(a) sketches the shape of $\ds\rho_{\eps}(x)=\int_V f(x,v)\ud v$. As $\eps\to0$, $\rho_{\eps}$ (solutions to kinetic system) approach to $\rho_0$, the solution of the limiting Keller-Segel system. Figure \ref{fig:conv_beforeBU}(b) shows that the convergence order of $f_{\eps}\to \rho_\eps F$ is almost first order (around $0.85$) in $l^2$ norm.

\begin{figure}
\centering
\subfigure[$\rho_\eps(x)$.]{ \includegraphics[width=0.45\textwidth]{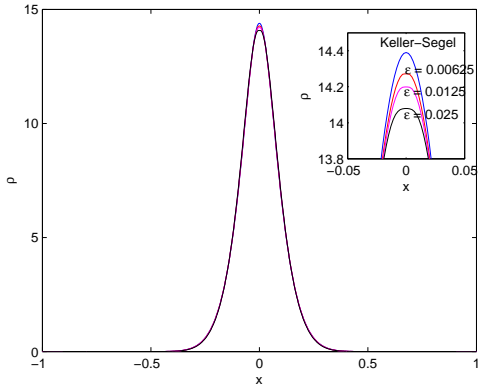}}
\subfigure[Convergence in $\eps$.]{ \includegraphics[width=0.45\textwidth]{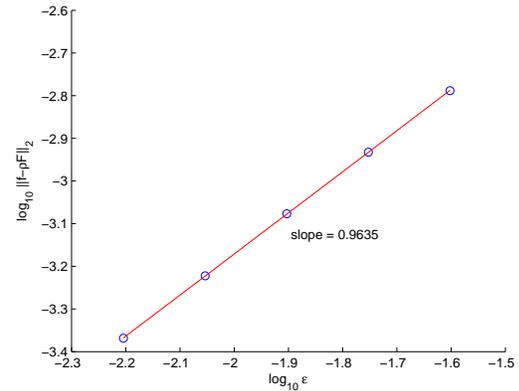}}
 \caption{{1D nonlocal model. The solutions of kinetic system and Keller-Segel system at time $t = 0.01$, before blow up time. The super-critical mass $M=4\pi$ is used. The left figure shows $\rho_{\eps}$ (for different $\eps$) and $\rho$. The right figure gives the convergence order of $||f_\eps(x,v) - \rho_\eps(x) F(v)||_2$.}}
\label{fig:conv_s}
\end{figure}

For a subcritical case, the initial data are taken to be (\ref{num:init_general}) with $M=\pi$. After a relatively long time $t=0.01$, we obtain the very
similar results. See figure \ref{fig:conv_s}. Now we have first order convergence in $\eps$, in $l^2$ norm.

\subsubsection{The stationary solution of kinetic system}

\begin{figure}
\centering
\includegraphics[width=0.8\textwidth]{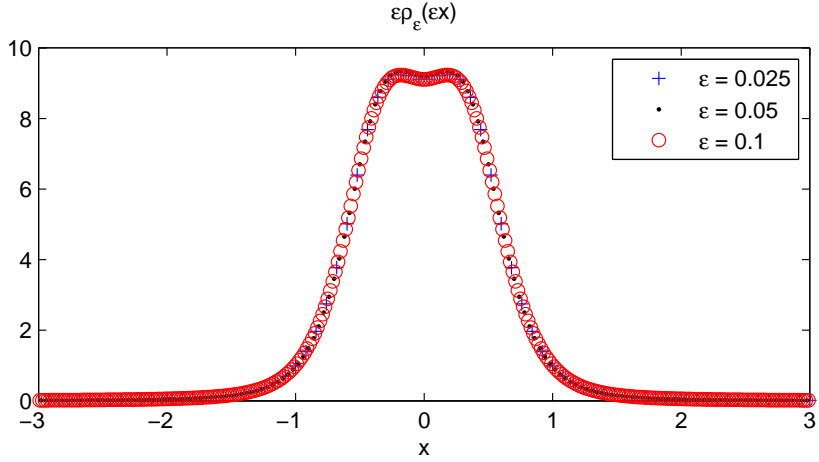}
\caption{{1D nonlocal model. The function $\eps\rho_{\eps}(\eps x)$
for different $\eps$. For supercritical mass $M=4\pi$ at time $t =
0.1\gg t_b$.}} \label{fig:erhoex}
\end{figure}

\begin{figure}
\centering
\includegraphics[width=0.8\textwidth]{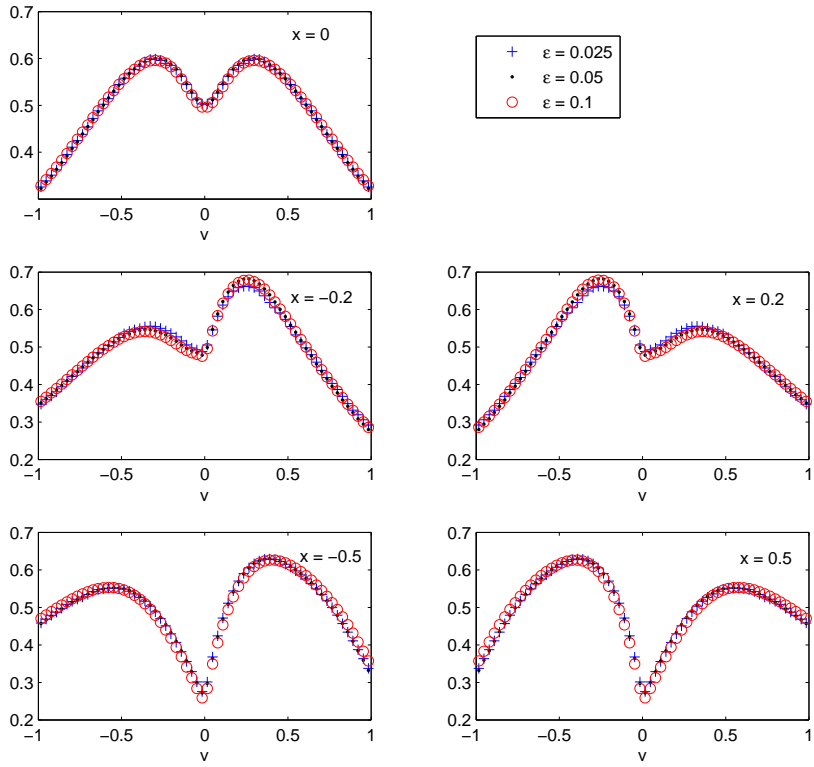}
\caption{{1D nonlocal model. The snapshots of the stationary state $\tilde F(x,v) =
\frac{f(\eps x,v)}{\rho(\eps x)}$ at certain locations $x$. For
supercritical mass $M=4\pi$ at time $t = 0.1\gg t_b$.}}
\label{fig:fonv}
\end{figure}

In this subsection we will study the solution of kinetic system
with large initial mass $M=4\pi$ at a relatively long time
$t=0.1$. After a long time, the solution stabilizes toward a
stationary state. This has not been proved nor intuitively
discussed in the literature. Figure \ref{fig:erhoex} shows the
function $\eps\rho_{\eps}(\eps x)$ at time $t=0.1$. Figure
\ref{fig:fonv} shows the function $\frac{f(\eps x,v)}{\rho(\eps
x)}$ as a function of $v$ at different $x$. These two figures
suggest the stationary state satisfies
\begin{equation}
\label{fmla: stnry} f_\infty^\epsilon(x,v) =
\frac{1}{\eps}\tilde\rho_\infty(\frac{x}{\eps})\tilde
F_\infty(\frac{x}{\eps},v).
\end{equation}
for some functions $\tilde\rho_\infty(x)$ and
$\tilde{F}_\infty(x,v)$. Therefore, let us consider the ansatz
\begin{eqnarray*}
\tilde{f}(t,x,v) &=& \eps f_{\eps}(t,\eps x,v) \\
\tilde{\rho}(t,x) &=& \eps \rho_{\eps}(t,\eps x) \\
\tilde{S}(t,x) &=& -\frac{1}{\pi} \log \abs{x} \ast \rho = S_{\eps} (\eps x) + C_{\eps}
\end{eqnarray*}
where $C_{\eps} = \frac{M}{\pi}\log{\eps}$, with $M$ the total
mass. The rescaled variables satisfy the following equations
\begin{equation*}
\eps^2 \frac{{\partial \tilde{f} }}{{\partial
t}} +  v  \frac{{\partial \tilde{f} }}{{\partial
x}}   = \left(F(v) +  \delta \tilde{S}(x,v) \right)\tilde{\rho} - \left(1+\int \delta \tilde{S}(x,v') \ud v'\right)\tilde{f}
\end{equation*}
\begin{equation*}
\tilde{S} = -\frac{1}{\pi}\log|x|\ast\tilde{\rho}
\end{equation*}
where $t>0$,
$x\in\Omega_{\eps}=[-\frac{x_{\max}}{\eps},\frac{x_{\max}}{\eps}]$,
$v\in V = [-v_{\max},v_{\max}]$. As $t\to \infty$,
$\frac{{\partial \tilde{f} }}{{\partial t}} \to 0 $, we have the
stationary solution $f_\infty^\epsilon$ should solve
\begin{equation}
\label{eq:ki_tilde_stati} v  \frac{{\partial f_\infty^\epsilon
}}{{\partial x}}   = \left(F(v) +  \delta \tilde{S}_\infty(x,v)
\right)\tilde{\rho}_\infty - \left(1+\int \delta
\tilde{S}_\infty(x,v') \ud v'\right)\tilde{f}
\end{equation}
where $\tilde{S}_\infty$ is obtained from $\tilde{\rho}_\infty$.
Clearly $\tilde F(x,v)$ satisfies
\[
\int_V \tilde F(x,v) \ud v = 1 \qquad \mbox{ and } \qquad \int_V
v\tilde F(x,v) \ud v = 0.
\]
The second identity can be derived by integrating
(\ref{eq:ki_tilde_stati}) over velocity space. These two
identities are also verified numerically. $\tilde \rho(x) = \eps
\rho(\eps x)$ is shown in Figure \ref{fig:erhoex}. Some snapshots of $\tilde F(x,v)$
are shown in Figure \ref{fig:fonv}.

\subsubsection{The interaction between several peaks}

It has been shown the interaction between several peaks for the
modified Keller-Segel system (\ref{eq:KellerSegel}) in
\cite{BCC08} by means of optimal transportation methods. Here we
will check this interaction for the kinetic system.

\

\begin{figure}
\centering
\includegraphics[width= 0.9 \textwidth]{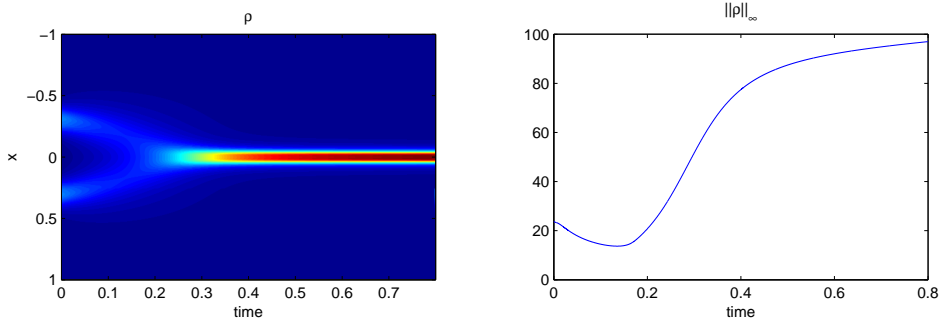}
\caption{{1D nonlocal model. The attraction between peaks. $\eps=0.1$. The total mass is $3\pi$. The critical mass is $2\pi$. $Nx = 400$. }}
\label{fig:2peaks_sym_3pi}
\end{figure}

\textbf{Case I: Two symmetric peaks, without enough mass in each
peak.-} The initial condition is taken as the sum of two
gaussian-like peaks,
\begin{equation}
\label{num:init_2peak}
f^I_{\eps}(x,v) = C\left(\frac{1}{2}e^{-80(x-0.3)^2}+\frac{1}{2}e^{-80(x+0.3)^2}\right) F(v)
\end{equation}
with a suitable constant $C$ such that the total mass is $3\pi$.
We take $\eps=0.1$. Figure \ref{fig:2peaks_sym_3pi} shows the time
evolution of $\rho_{\eps}$ (left) and $||\rho_{\eps}||_{\infty}$
(right). The solution starts with diffusion separately, and then
the peaks merge. At beginning there is no enough mass in either
peak to concentrate ($\frac{3\pi}{2}<M_{c}=2\pi$). But after
they merge, the total mass in the new peak is large enough to form
an aggregation.

\

\textbf{Case II: Two symmetric peaks, with enough mass in each
peak.-} We take the same initial setting (\ref{num:init_2peak}),
with a different constant $C$ such that the total mass is $5\pi$.
Now the mass in each peak is large enough to concentrate
($\frac{5\pi}{2}>M_{c}=2\pi$).
\begin{figure}
\centering
 \includegraphics[width=0.9\textwidth]{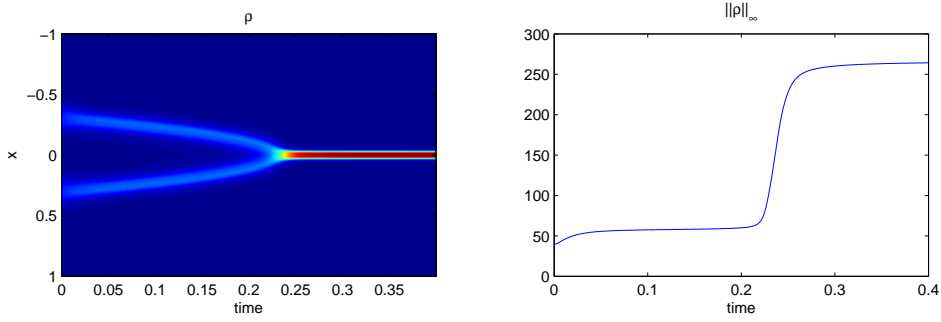}
\caption{{1D nonlocal model. The attraction between peaks. $\eps=0.05$. The total mass is $5\pi$. The critical mass is $2\pi$. $Nx = 400$. }}
\label{fig:2peaks_sym_5pi}
\end{figure}
In Figure \ref{fig:2peaks_sym_5pi}, we can see two distinguished
different phases during the time evolution. The first phase
consists in the appearance of a metastable state where the
concentration in each peak is formed but slowly moving toward each
other. Then the two peaks merge, which increases the total mass
into the final larger peak. It continues to aggregate the mass
around it and finally reaches the stationary state.

\

\begin{figure}
\centering
  \includegraphics[width=0.9\textwidth]{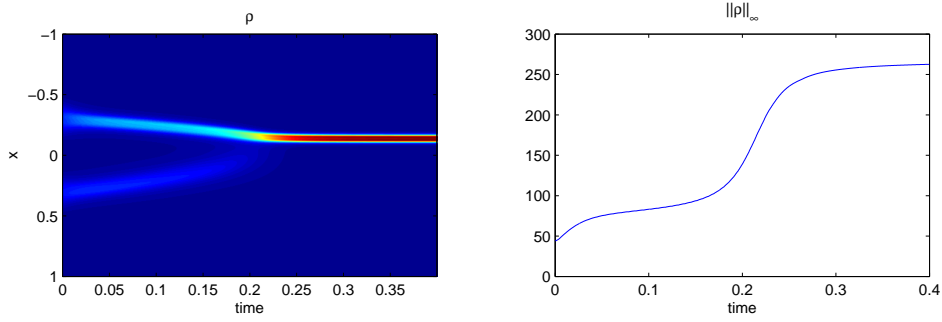}
\caption{{1D nonlocal model. The attraction between peaks. $\eps=0.05$. The total mass is $5\pi$ with $2.2\pi$ near $x = 0.3$, $2.8\pi$ near $x=-0.3$. The critical mass is $2\pi$. $Nx = 400$. }}
\label{fig:2peaks_unsym_5pi}
\end{figure}

\textbf{Case III: Two asymmetric peaks, with enough mass in each
peak.-} We take a nonsymmetric initial setting
\begin{equation*}
f^I_{\eps}(x,v) = C\left(2.2e^{-80(x-0.3)^2}+2.8e^{-80(x+0.3)^2}\right)F(v)
\end{equation*}
with a suitable constant $C$ such that the total mass is $5\pi$.
The mass in each peak is still large enough to concentrate and we
see that the mass tends to concentrate around the center of mass
located this time much closer to the first peak due to asymmetry
in Figure \ref{fig:2peaks_unsym_5pi}.

\

\begin{figure}
\centering
 \includegraphics[width=0.9\textwidth]{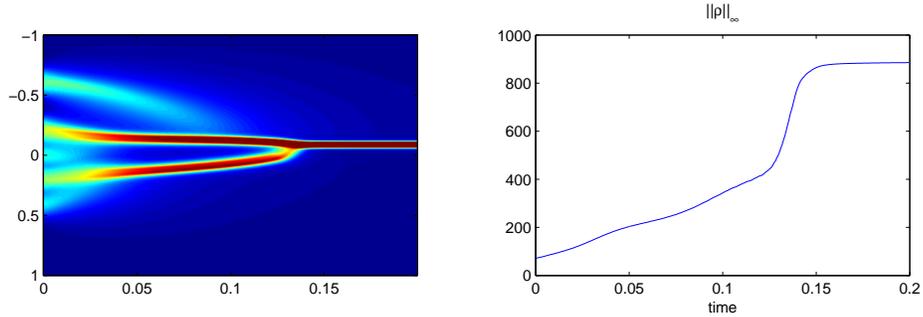}
\caption{{1D nonlocal model. The interaction between five peaks. $\eps=0.05$. The total mass is $11\pi$. The critical mass is $2\pi$. $Nx = 400$. }}
\label{fig:5peaks_unsym_11pi}
\end{figure}

\textbf{Case IV: Five unsymmetric peaks, with total mass
$M>5M_{c}$.-} We take a nonsymmetric initial setting, with five
different peaks,
\begin{equation*}
\ds f^I_{\eps}(x,v) = C\sum_{j=1}^{5} w_j e^{-160(x-c_j)^2}
\end{equation*}
with $w = [0.5, 1.2, 0.8, 0.6, 1]$, $c = [0.4, -0.2, -0.6, 0,
0.2]$. Here we pick up a suitable constant $C$ such that the total
mass is $11\pi$. Now, we observe in Figure
\ref{fig:5peaks_unsym_11pi} the complicated dynamics of merging
between the different peaks before forming the final peak at the
center of mass and converging toward the stationary state.

\subsection{The 1D local model: blow up in finite time}

We numerically check the open problem of the blow up property of
the solution to the one dimensional kinetic local model
(\ref{eq:BC1d}). As mentioned before, a theoretical prediction on
this blow up is still lacking, see \cite{S11}. We consider the
initial data given by (\ref{num:init_general}). The critical mass for corresponding Keller-Segel model is again
$$M_{c} = 2\pi.$$

For the super-critical case, we take total mass $M = 5\pi$.

\subsubsection{The problem in the convergence with fixed grids}

In this test $\eps = 0.4$.  Different grid sizes in $x$ are used.
The number of grid points in each of the simulations are $N_x=
250,500,1000,2000,4000,8000$ respectively and the time step sizes
are taken as $\Delta t = \eps \Delta x/v_{\max}$. We still
consider the convergence order defined as in (\ref{num:err}).

\begin{figure}
\centering
  \includegraphics[width=0.8\textwidth]{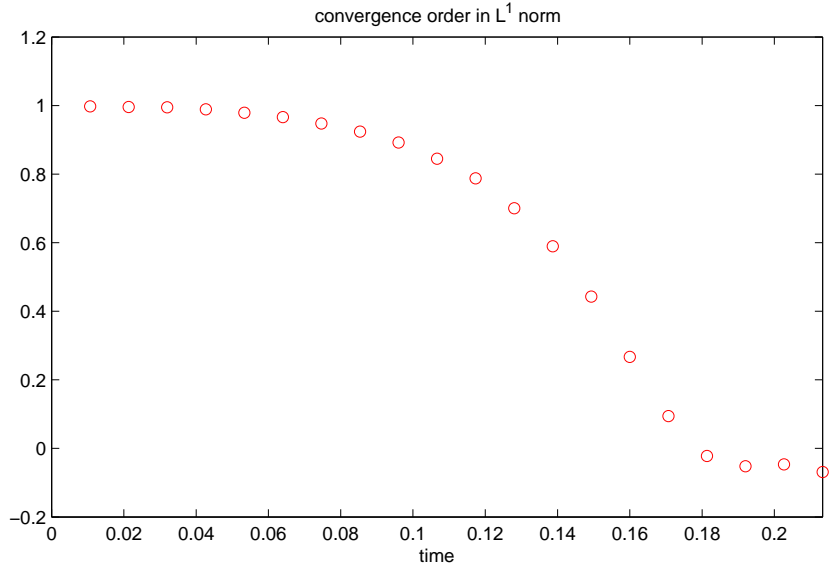}
\caption{{1D local model. The convergence order of the solutions to
(\ref{eq:BC1d}) at different time. $\eps = 0.4$. The total mass is
$5\pi>M_{c}$.  }} \label{fig:convorder_fix}
\end{figure}

Figure \ref{fig:convorder_fix} shows the convergence order in
$L^1$ norm at different times $\tau_k = k t_{\max}/20$, where
$1\le k\le 20$. The solutions show a first order convergence when
$t<0.12$, then the convergence order decreases as time evolves.
After $t>0.16$ a negative convergence order is seen, which means
the scheme is not convergent after that time. One cannot improve
these convergence orders by using a finer grid. Figure
\ref{fig:convorder_fix} strongly suggests that the blowup happens
in the solution during $0.12<t<0.16$. Next we investigate this
time period by using the adaptive grids proposed in section
\ref{sec:adaptivegrid}.

\subsubsection{The convergence with adaptive grids}

\begin{figure}
\centering
  \includegraphics[width=0.8\textwidth]{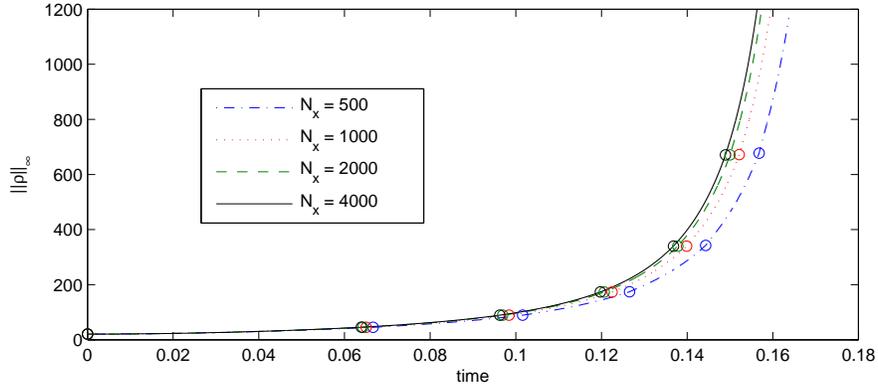}
\caption{{1D local model. The time evolution of $||\rho||_{\infty}$ by
using adaptive grids, with different initial grids. The circles
shows the time steps when the grids are doubled. $\eps = 0.4$. The
total mass is $5\pi>M_{c}$.  }} \label{fig:conv_adap}
\end{figure}

Now we use the adaptive grids proposed in section
\ref{sec:adaptivegrid}. We start the simulation with different
grid sizes $N_x = 500, 1000, 2000, 4000$. Figure
\ref{fig:conv_adap} shows the time evolution of
$||\rho||_{\infty}$. A nice convergence toward some density
blow-up is observed.

\subsubsection{The convergence as $\eps\to0$}

\begin{figure}
\centering
  \includegraphics[width=0.8\textwidth]{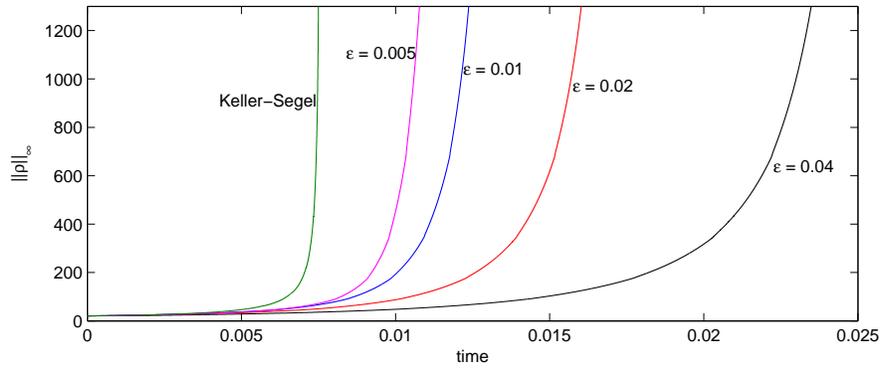}
\caption{{1D local model. The time evolution of $||\rho||_{\infty}$ by using adaptive grids, with different $\eps$. The total mass is $5\pi>M_{c}$. Initially $N_x = 1000$. }}
\label{fig:conv_adap_eps}
\end{figure}

Next we study the convergence as $\eps\to0$. We apply the adaptive
grids described above and take $N_x=1000$ at beginning. We take
total mass to be $M = 5\pi > M_{c}$. Since the solutions to
kinetic equations also blow up in finite time, they converge to
the solution of Keller-Segel system in a totally different way.
Figure \ref{fig:conv_adap_eps} shows the time evolution of
$||\rho||_{\infty}$ of kinetic equations with different $\eps$, as
well as that Keller-Segel system. The way of asymptotic
convergence is totally different with Figure \ref{fig:kin_KS},
where the solutions of kinetic equations are globally bounded.

\subsubsection{Subcritical case: long time behavior}

Now we check the long time behavior of subcritical case. Let us
remark that in the case of the limiting Keller-Segel model, it is
known that the long time asymptotics should be given by a
self-similar solution whose profile is dictated by a stationary
scaled problem, see \cite{BCC08} for instance. We now check this
point for the 1D local kinetic model by taking as total mass
$M=\pi$. The simulation is performed on $x\in[-10,10]$. A fixed
grid with $N_x = 2000$ is used, with $\eps = 0.2$. We consider the
same change of variables as for the Keller-Segel model,
$$
\ds R(t) = \sqrt{1+2t}, \quad y = \frac{x}{R(t)}, \quad \tau = \log R(t),
$$
\begin{equation}
\label{def:rhotilde_BCsub}
\ds\rho(t,x)=\frac{1}{R(t)} \tilde{\rho}(\tau,y).
\end{equation}
Here $\rho$ is the density we computed from numerical simulation.

\begin{figure}
\centering
 \subfigure[$\tilde{\rho}(\tau,y)$.]{  \includegraphics[width=0.65\textwidth]{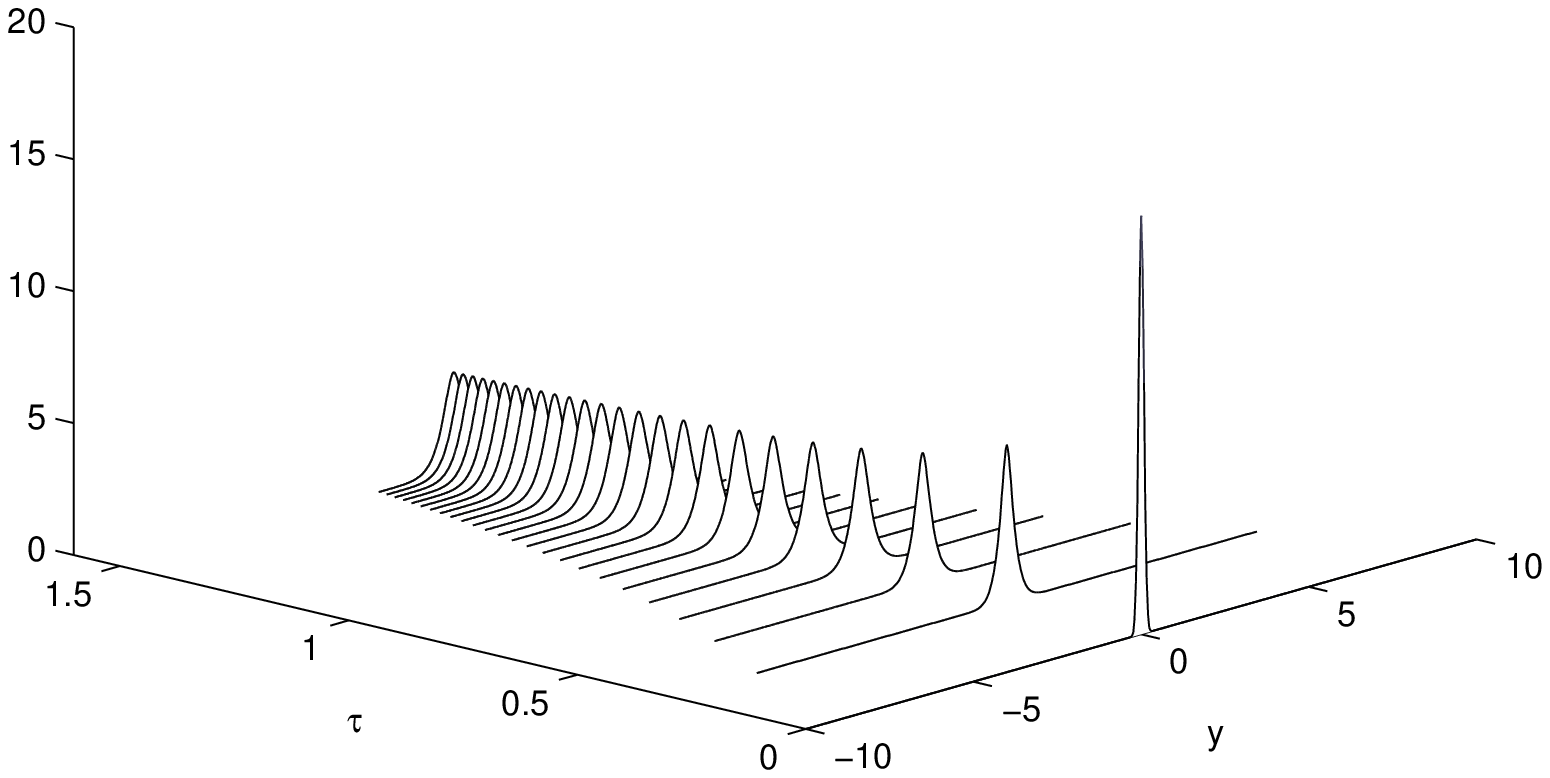}}
 \subfigure[$||\tilde{\rho}(\tau,y)-\tilde{\rho}_{\infty}||_1$.]{  \includegraphics[width=0.32\textwidth]{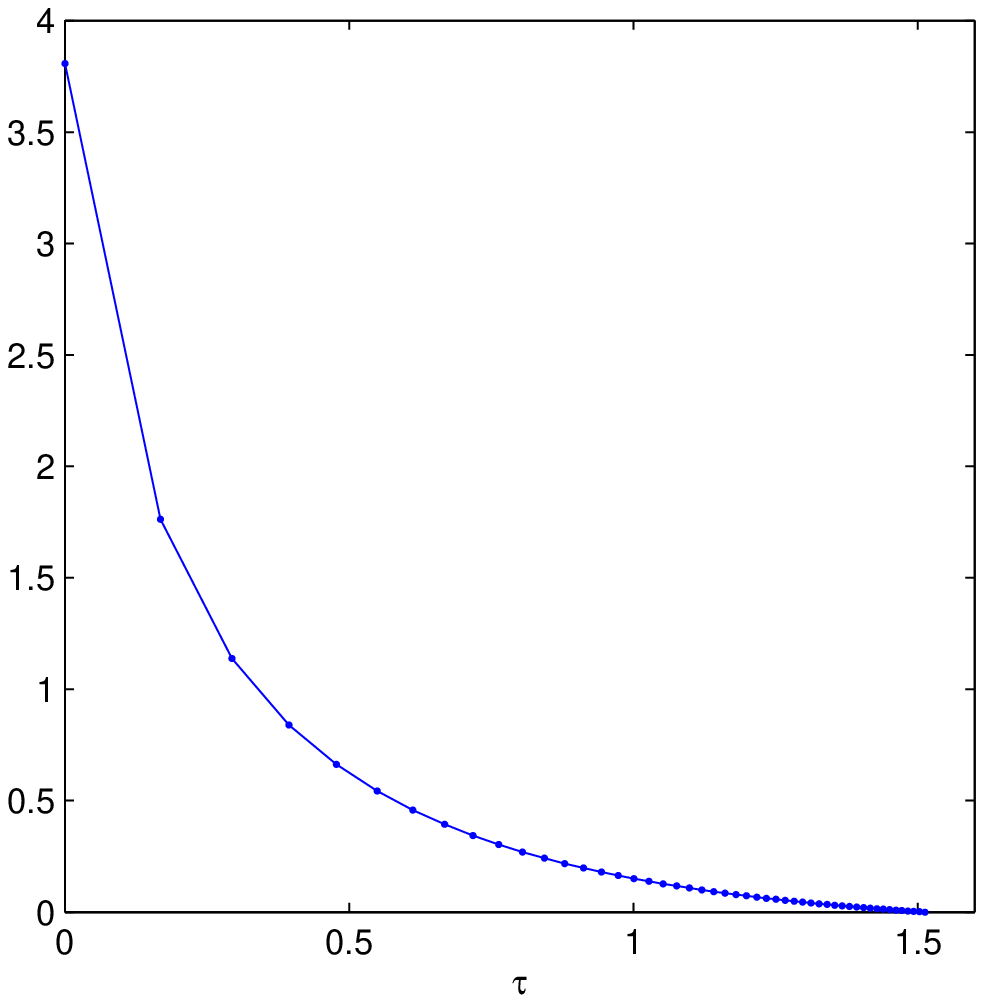}}
\caption{{1D local model. The time evolution of $\tilde\rho$ for
the subcritical mass $M=\pi<M_{c}$. $N_x = 2000$. $\eps = 0.2$.
}} \label{fig:BC1d_sub}
\end{figure}

Figure \ref{fig:BC1d_sub}(a) shows the time ($\tau$) evolution of
$\tilde{\rho}$ until time $\tau=1.5$ ($t=10$) as a function of $y$. It
strongly suggests that $\tilde{\rho}$ is approaching some
stationary state. We denote it by $\tilde{\rho}_\infty$. In Figure
\ref{fig:BC1d_sub}(b) we show the time evolution of the relative
error $||\tilde{\rho}(\tau,y)-\tilde{\rho}_\infty||_1$, with
$\tilde{\rho}_\infty$ approximated through
(\ref{def:rhotilde_BCsub}) at time $\tau=10$. Clearly the
solutions are approaching the stationary state, in other words,
\[\ds \Big|\Big|\rho(t,x) - \frac{1}{R(t)} \tilde{\rho}_\infty\left(\frac{x}{R(t)}\right) \Big|\Big|_1\to0, \quad \mbox{ as } t\to\infty.\]

\subsection{The 2D local model: solution behavior between theoretical thresholds}

In the final simulation we test the 2D local model. As shown in
\cite{BC09}, the solution blows up with total mass larger than the
critical mass
$$
M_c = \frac{32\pi}{|V|} = 32,
$$
where $|V| = \pi \omega_{\max}^2 = \pi$. And the global solution exists with the total mass lower than the other threshold,
$$
m_c = \frac{0.806\pi}{|V|} = 0.806.
$$
Besides, as mentioned in the introduction, the critical mass for the corresponding Keller-Segel model (\ref{eq:KellerSegel_tilde}) is
$$
M_{KS} = \frac{8\pi D}{\chi} = 16.
$$
The gap between the two estimates $M_c$ and $m_c$ should be much
smaller as observed from our numerical simulation. It is difficult
to ascertain based on numerical simulations if the kinetic system
show a clear dichotomy as in the Keller-Segel system.

We take $r_{\max} = 2$, $\omega_{\max} = 1$, with the equilibrium given by
$$
\ds F(\omega) = \frac{1}{2\pi \int_0^{\omega_{\max}}\omega \ud \omega} = \frac{1}{\pi}.
$$
The grid points in each directions are $N_r = 1000$, $N_\omega = 32$, $N_\theta = 32$. The
initial data are taken as
$$
\ds \tilde\rho^I (r) = C r e^{-15r^2}, \quad h^I(r,\omega,\theta) = \tilde \rho^I(r) F(\omega),
$$
where $C = C(M)$ is a constant determined by the total mass $M$.

\begin{figure}
\centering
 \includegraphics[width=0.95\textwidth]{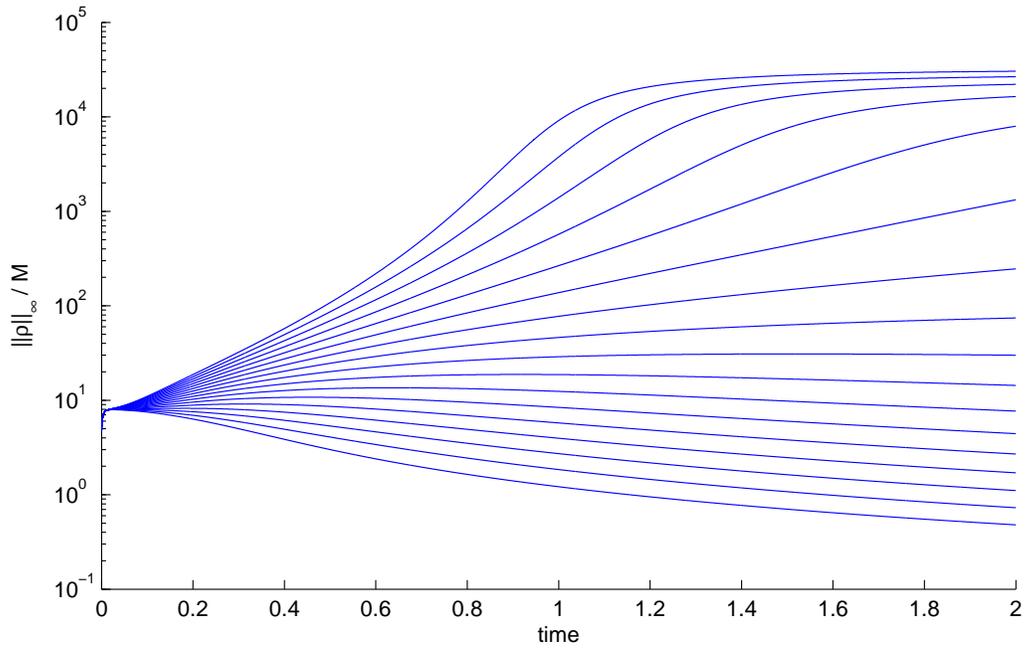}
\caption{{2D local model. The time evolution of $\ds\frac{||\rho||_\infty}{M}$ for the 2D local model with total masses M. From up to down: $M = 33, 31, 29,\dots,3,1$. $N_r = 1000$. $\eps = 1$. }}
\label{fig:BC2d_maxrho}
\end{figure}

Figure \ref{fig:BC2d_maxrho} shows the time evolution of
$||\rho||_\infty/M$, for different total masses $M$. We take $\eps
= 1$. It suggests that the global solutions exists for $M\le17$,
which is much bigger than the theoretical thresholds $m_c$. While
for mass $M\ge25$, the solution blows up, even if the total mass
belongs to the range of masses in which there are no theoretical
results \cite{BC09}. Let us comment that $||\rho||_\infty$ has a
upper bound due to the limitation of grid size. When we use a
finer grid, this upper bound would be larger.

\begin{figure}
\centering
 \includegraphics[width=0.95\textwidth]{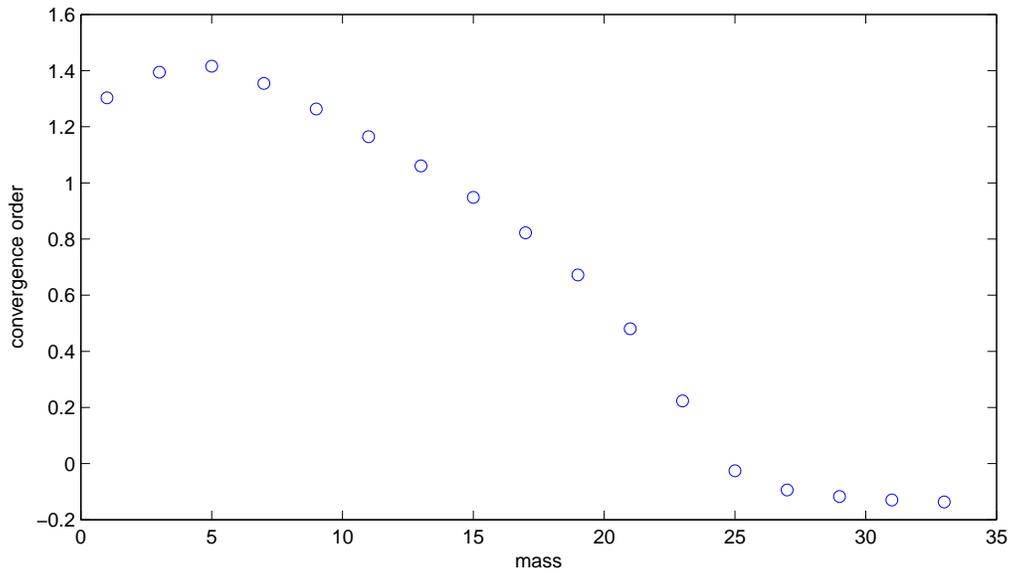}
\caption{{2D local model. The convergence order (in $L^1$ norm) of the solutions to the 2D kinetic system at time $t=2$ as $\Delta r\to0$, with different total masses $M$.  $\eps = 1$. }}
\label{fig:BC2d_order}
\end{figure}

Finally we compute the convergence of the solutions to the kinetic
system as $\Delta r\to0$, with different total masses $M$. As
shown in Figure \ref{fig:BC2d_order}, the numerical scheme shows a
fist order convergence for $M\le15$ (note that the critical mass
of the Keller-Segel model is $M_{KS} = 16$). While for $M\ge25$,
the numerical solutions do not convergence, which suggests that
the solutions blow up.

\section*{Acknowledgements}

This work was partially supported by NSF grant No. DMS-0608720, NSF FRG grant DMS-0757285. JAC was partially supported by the projects  Ministerio de Ciencia
e Innovaci\'on MTM2011-27739-C04-02 and 2009-SGR-345 from
Ag\`encia de Gesti\'o d'Ajuts Universitaris i de
Recerca-Generalitat de Catalunya. BY would like to thank Professor
Shi Jin in University of Wisconsin-Madison for the guidance and
fruitful discussion. BY also would like to thank the warm
hospitality at the Department of Mathematics in Universitat
Aut\`{o}noma de Barcelona during his visit. Both authors would
like to thank IPAM where part of this work was initiated during
the kinetic theory program.

\bibliographystyle
{plain}
\bibliography{chemoref}
\end{document}